\documentclass[10pt,a4paper]{article}
\usepackage{amsmath}
\usepackage{amscd}
\usepackage{amssymb}

\newtheorem{thm}{Theorem}
\newtheorem{defn}{Definition}
\newtheorem{cor}{Corollary}
\newtheorem{lem}{Lemma}
\usepackage{color}
\usepackage{colordvi}

\newtheorem{prop}{Proposition}
\newtheorem{ex}{Example}
\newtheorem{rk}{Remark}

\usepackage{amssymb}

\newcommand{\C}{\mathbb C}
\def \Ci {{C^\infty}}
\def \Cl {{C\ell}}
\def\cutoffint{-\hskip -10pt\int}

\def\dbar{d{\hskip-1pt\bar{}}\hskip1pt}
\newcommand{\e}{\varepsilon}

\def \endsquare{ $\sqcup \!\!\!\! \sqcap$ }

\newcommand{\N}{\mathbb N}
\def\otherterm#1{{\it#1}}

\def\pdo{$\Psi${\rm DO}$\quad$}
\def\pdos{$\Psi${\rm DOs}}

\newcommand{\R}{\mathbb R}

\newcommand{\Z}{\mathbb Z}

\begin{document}
\title{ \bf
The noncommutative residue and  canonical trace in
  the light of Stokes'  and continuity properties}
\author{  Sylvie Paycha  }
\maketitle
\section*{Abstract}  We  show that  the noncommutative
residue density, resp.  the cut-off regularised integral are the only  {\it
  closed} linear,
resp.  continuous closed linear  forms on certain classes of symbols. This leads to  alternative proofs of the   uniqueness of the
noncommutative residue, resp. the canonical trace as  linear, resp. continuous
linear forms  on certain classes of
classical pseudodifferential operators which vanish on brackets.\\ The
uniqueness of the canonical trace actually holds  on classes of  classical pseudodifferential with
{\it vanishing residue density}  which include  non integer order
operators in all dimensions and odd-class (resp. even-class)  operators in
odd (resp. even) dimensions. The description   of the canonical trace for non
integer order operators as an
 integrated global density on the manifold   is
 extended to  
odd-class (resp. even-class) operators in odd (resp. even) dimensions on the
grounds of defect formulae for regularised traces of classical
pseudodifferential operators. 
\section*{Acknowledgements}
Let me thank Matthias Lesch for stimulating discussions around the canonical
trace as well as Elmar Schrohe for his precious comments on these uniqueness
issues. I am also  endebted to Maxim Braverman whom I thank  for crucial remarks on an earlier
version of this paper.  I further thank Mikhail Shubin for his constant
 encouragements and enlightening comments on trace related issues.  Last but
 not least, I am very grateful to  Carolina Neira and  Marie-Fran\c coise Ouedraogo 
 for a careful reading of a former version of these notes which greatly helped me
 improve the presentation.
 \section*{Introduction}
  The uniqueness of the
noncommutative residue  originally introduced by Adler and Manin in the one
dimensional case was then extended to all dimensions by Wodzicki in \cite{W1} (see also \cite{W2} and \cite{K}
for a review) and
proved independently by Guillemin \cite{G2}. Since then
  other proofs, in particular a homological proof on symbols
in \cite{BG}   and various extensions
of this uniqueness result were derived, see \cite{FGLS} for a
generalisation to manifolds with boundary, see \cite{S} for a generalisation to
manifolds with conical singularities (both of which prove uniqueness up to
smoothing operators),  see \cite{L} for an extension to
log-polyhomogeneous operators as well as for an argument due to Wodzicki to
 get uniqueness on the whole algebra of classical operators, see \cite{Po2}
 for an extension to Heisenberg manifolds. \\ In contrast to the  familiar
characterisation of the noncommutative residue as the unique trace on the
algebra of all classical pseudodifferential operators, only
recently was the focus \cite{MSS} drawn on  a characterisation of the
canonical trace
as  the unique linear extension of the
ordinary trace to non integer order classical pseudodifferential operators
which vanishes on non integer order brackets\footnote{The authors of
  \cite{MSS} actually extended the uniqueness to odd-class,
  resp. even-class operators in odd, resp. even dimensions.}. \\
We  revisit  and
slightly improve these results handling  the noncommutative residue and the  canonical trace on an equal
footing via a characterisation of closed linear forms
on certain classes of symbols. 
 \\ \\
  A cornerstone in
 our approach  is the requirement that a
linear form  satisfies  {\it  Stokes' property } (or
be {\it closed} in the language of noncommutative geometry) on a certain class of symbols  i.e. that it vanishes on partial
derivatives in that class. The vanishing on derivatives is  a  natural requirement in  view of the fact that
any distribution on $\R^n$ with vanishing derivatives is proportional to the ordinary
integral on $\R^n$; it serves here to characterise its unique closed extension
given by the cut-off regularised integral $\cutoffint_{\R^n}$ defined by Hadamard
finite parts. This leads to a  characterisation of the noncommutative
residue on symbols (Theorem \ref{thm:uniquenessres}) on the one hand and the  cut-off regularised
integral  on symbols  (Theorem \ref{thm:uniquenesscutoffint}) 
on the other hand. \\ The link between the vanishing on brackets of a linear
functional on classical operators and the vanishing on partial derivatives of
a linear functional on symbols can best be seen from the simple formula $[x_i,
{\rm Op}(p)]= -i\, {\rm Op}\left(\partial_{\xi_i} p\right)$ for any symbol $p$,
which lies at the heart of the proof of the uniqueness of the canonical trace
in \cite{MSS}. We  deduce  from there the uniqueness of a linear form on
classical
pseudodifferential operators  which
vanishes on brackets from  the  uniqueness of a linear form on classical
pseudodifferential symbols which vanishes on partial derivatives in the
Fourier variable $\xi$, leading to   characterisations of  the residue
(Theorem \ref{thm:uniquenessresM}) on the one hand and the
canonical trace  (Theorem \ref{thm:uniquenesstraceM}) on the other hand. \\ To characterise such functionals we heavily rely on
the fact that any homogeneous symbol with  vanishing residue   can be written
as a sum of derivatives \cite{FGLS} and coincides up to a smoothing symbol
with a sum of derivatives of symbols whose order is $1+$ the order of the
original  symbol. This is why we then  consider classes of
operators with {\it vanishing residue density}  in order to carry out the linear
extensions. \\ \\
 The  vanishing of the residue density which therefore plays a crucial
 role for uniqueness issues,  arises once more for existence
 issues. This  indeed turns out to be an essential ingredient
 in section 2, where we show that the canonical trace is well defined as an
 integrated global density  on certain classes of classical
 pseudodifferential operators, such as odd-class operators in odd dimensions
 and even-class operators in even dimensions. To do so, we approximate the operator
under study  along a holomorphic path of classical operators and use a defect
formula for regularised traces derived in \cite{PS}. \\ This is carried out
along  of a  line of thought underlying  Guillemin's \cite{G2},  Wodzicki's
\cite{W2} and later Kontsevich and Vishik's \cite{KV} work (see also
\cite{CM}); a classical \pdo $A$ is embedded in a holomorphic
family $z\mapsto A(z)$ with $A(0)=A$,   the
canonical trace  of which yields a meromorphic map $z\mapsto {\rm
  TR}(A(z))$. These authors focus on the important case of $\zeta$-regularisation
$A^Q(z)=A\, Q^{-z}$ built from some admissible elliptic operator $Q$   with positive
order. In particular,  if $A$ has  non integer order  then   $z\mapsto {\rm TR}(A^Q(z))$ is holomorphic
at $z=0$, the canonical trace density of $A$ is globally defined and integrates
over $M$ to the canonical trace TR$(A)$  of $A$ which  coincides with 
$ \lim_{z\to 0}{\rm TR}(A^Q(z))$  independently  of the choice of
$Q$. \\ Similar continuity results hold for odd-class (resp. even-class)
operators $A$ in odd (resp. even)
dimensions; it was  observed in \cite{KV}  (resp. \cite{Gr}) that for an odd-class elliptic operator $Q$   with even positive 
order close enough to a positive
self-adjoint operator, and $A$ and odd-(resp. even-) class operator in odd
(resp. even) dimensions, 
 the map $z\mapsto {\rm TR}(A^Q(z))$ is holomorphic at $z=0$ and
${\rm Tr}_{(-1)}(A):= \lim_{z\to 0}{\rm TR}(A^Q(z))$ is independent of the choice
of $Q$. \\ As a straightforward application of  defect formulae both on the symbol
and the operator level derived in
\cite{PS}, we extend these results to {\it any holomorphic  family $A(z)$} with non
constant affine order  such that $A=A(0)$ and  $A^\prime(0)$ lie
in the odd- (resp. even-) class. We infer from there that in odd (resp. even) dimensions
\begin{enumerate}
\item  the map  $z\mapsto {\rm TR}(A(z))$ is  holomorphic at $z=0$, 
\item  the canonical trace density  is globally defined for  any  odd-
  (resp. even-) class operator $A$, and
  integrates over $M$ to the canonical trace TR$(A)$, 
\item 
${\rm TR}(A)=  \lim_{z\to 0} {\rm TR}(A(z))$ 
independently of any appropriate (see above initial conditions) choice of
the family  the family $A(z)$.
\end{enumerate} This  shows in particular that both Kontsevich and Vishik's
(resp. Grubb's) extended trace
${\rm Tr}_{(-1)}$ on odd- (resp. even-) class operators  in odd (resp. even) dimensions  and  the
symmetrised trace ${\rm Tr}^{\rm sym}$  introduced by Braverman
in \cite{B}  on odd-class operators in odd dimensions coincide with the canonical trace TR. \\\\ 
To sum up, the characterisation we provide of  the noncommutative residue and of
the canonical trace  on the grounds of  a characterisation of closed linear forms on
certain classes of symbols sheds
light on common mechanisms that underly their uniqueness. It  brings out
the importance of the  closedness requirement on the underlying functionals on
symbols, which was already  implicit in the homological proofs of the
uniqueness of the residue. In the case of the canonical trace it  further puts forward the role of the vanishing
of the residue on the   symbol level and of the residue density on the
 operator level which also turns out to be  an essential ingredient for  existence issues. 
\\ \\ The paper is organised as follows:
\begin{enumerate}\item {\bf  Uniqueness: characterisation  of  closed linear forms on symbols}
\begin{enumerate} \item Notations
\item Stokes' property versus translation invariance
\item A characterisation of the noncommutative residue and its kernel
\item  A characterisation of the cut-off regularised integral
  $\cutoffint_{\R^n}$ in terms of Stokes' property
\end{enumerate}

\item {\bf Existence: The canonical trace on 
odd- (resp. even-) class  operators in odd (resp. even) dimensions}
\begin{enumerate}
\item Notations
\item Classical symbol valued forms on an open subset
\item The noncommutative residue on classical pseudodifferential operators
\item The canonical trace on non integer order  operators
\item Holomorphic families of classical pseudodifferential operators
\item Continuity of the canonical trace on non integer order pseudodifferential
  operators
\item Odd-class (resp. even-class) operators embedded in holomorphic families
\item The canonical trace on  odd- (resp. even-)class operators in odd
  (resp. even) dimensions
\end{enumerate}
\item {\bf Uniqueness: Characterisation of  linear forms on operators that vanish on  brackets}
\begin{enumerate}
\item Uniqueness of  the noncommutative residue 
\item Uniqueness of the canonical trace
\end{enumerate}
\end{enumerate}
\vfill \eject \noindent
  \section{Uniqueness: Characterisation  of closed  linear forms on symbols }
We  characterise the noncommutative residue and the
cut-off regularised integral in terms of a closedness requirement on linear
 forms  on classes of classical symbols with constant coefficients
on $\R^n$.   \subsection{Notations} We only give a few definitions and refer the reader to
\cite{Sh,T,Tr} for further details on classical pseudodifferential
operators.\\  For any complex number $a$, let us denote by ${\cal S}^a_{\rm
  c.c}(\R^n)$ 
 the set of smooth functions on $\R^n$ called symbols with constant
 coefficients, such that for any multiindex
$\beta\in \N^n$ there is a constant $C(\beta)$ satisfying the following requirement:
$$\vert\partial_\xi^\beta \sigma(\xi)\vert\leq C(\beta) \vert (1+\vert
\xi\vert)^{{\rm Re}(a)-\vert \beta\vert}$$ where Re$(a)$ stands for the real
part of $a$, $\vert
\xi\vert$ for the euclidean norm of $\xi$. We single out the subset  $CS^a_{c.c}(\R^n)\subset {\cal
  S}_{\rm c.c}^a(\R^n) $ of   symbols $\sigma$, called classical symbols of order $a$  with constant
coefficients,   such that 
\begin{equation}\label{eq:asymptsymb}
\sigma(\xi)=
\sum_{j=0}^{N-1} \chi(\xi)\, \sigma_{a-j}( \xi) +\sigma_{(N)}( \xi)
\end{equation}
where $\sigma_{(N)}\in {\cal S}_{\rm c.c}^{a-N}(\R^n)$ and
 where $\chi$ is a smooth cut-off function which vanishes in a
  small
   ball of $\R^n$ centered at $0$ and which is constant equal to $1$ outside
   the unit ball. Here $\sigma_{a-j,}, j\in \N_0$ are positively homogeneous of degree
   $a-j$. \\The ordinary product of functions sends $CS^a_{c.c}(\R^n)\times
CS^b_{\rm c.c}(\R^n)$ to $CS^{a+b}_{c.c}(\R^n)$ provided $b-a\in \Z$; let   $$ CS_{c.c}(\R^n)= \langle
\bigcup_{a\in \C}CS^a_{c.c}(\R^n)\rangle $$ denote the algebra generated by   all classical
symbols with constant coefficients  on $\R^n$.  Let  $$CS_{c.c}^{-\infty}(\R^n)=  \bigcap_{a\in
  \C}CS^a_{c.c}(\R^n)$$ be the algebra of smoothing symbols. We write
$\sigma\sim \sigma^\prime$ for two symbols $\sigma, \sigma^\prime$  which differ by a smoothing symbol.\\
 We also
denote by $CS_{\rm
  c.c}^{<p}(\R^n):= \bigcup_{{\rm Re}(a)<p} CS_{\rm c.c}^a(\R^n)$, the set of classical
symbols of order with real part $<p$ and by   $CS_{\rm
  c.c}^{\notin \Z}(\R^n):= \bigcup_{a\in\C- \Z} CS_{\rm c.c}^a(\R^n)$
the set of non integer order symbols.
\\ We equip the set $CS^{a}_{\rm c.c}(\R^n)$ of classical symbols of order $a$ with a Fr\'echet structure with the help of the following
semi-norms labelled by multiindices $\beta$ and integers $j\geq
0$, $N$  (see \cite{H}):
\begin{eqnarray*}
&{} & {\rm sup}_{ \xi \in \R^n} (1+\vert \xi\vert)^{-{\rm Re}(a)+\vert \beta\vert} \, \Vert  \partial_\xi^\beta \sigma( \xi)\Vert;\\
&{}&  {\rm sup}_{ \xi\in \R^n}  (1+\vert \xi\vert)^{-{\rm Re}(a)+N+\vert \beta\vert}\Vert  \partial_\xi^{\beta}\left(\sigma-\sum_{j=0}^{N-1} \chi(\xi)\, \sigma_{a-j}\right)(\xi) \Vert;\\
&{}& {\rm sup}_{ \vert\xi\vert=1}  \Vert \partial_\xi^{\beta} \sigma_{a-j}( \xi) \Vert.
\end{eqnarray*}
$CS_{\rm c.c}^{-\infty}(\R^n)$ is equipped with the natural induced topology
so that a  linear $\rho$ which extends to continuous linear maps $\rho_a$  on
$CS_{\rm c.c}^a(\R^n)$ for any $a\in \Z\cap
]-\infty, -K]$ (with  $K$ some arbitrary positive number) is continuous.\\ \\
We borrow from \cite{MMP} (see also \cite{LP}) the notion of \pdo-valued form.
\begin{defn}Let $k$ be a non negative integer,
$a$ a complex number. We let
  \begin{eqnarray*}
\Omega^k\,CS^a_{\rm c.c}(\R^n)&= &\{\alpha \in \Omega^k(\R^n),\quad
 \alpha = \sum_{I\subset \{1, \cdots, n\}, \vert I
\vert = k} \alpha_{I} (\xi) \,  d\xi_I\\
 {\rm with } &{}& \quad \alpha_{I} \in CS_{\rm c.c}^{a-\vert I\vert}(\R^n)\}
\end{eqnarray*}
denote the set of order $a$ classical symbol valued forms on $\R^n$ with
constant coefficients.
Let  \begin{eqnarray*}
\Omega^k\,CS_{\rm c.c}(\R^n)&= &\{\alpha \in \Omega^k(\R^n), \quad
 \alpha = \sum_{I\subset \{1, \cdots, n\}, \vert J
\vert = k} \alpha_{I} (\xi) \,  d\xi_I\\
 {\rm with } &{}& \quad \alpha_{I} \in CS_{\rm c.c}(\R^n)\}
\end{eqnarray*} denote the set of classical
symbol valued $k$-forms on $\R^n$ of all orders with constant coefficients.
  \end{defn}
The exterior product on forms induces a
product  $\Omega^k CS_{\rm c.c}(\R^n)\times \Omega^l CS_{\rm c.c}(\R^n)\to
\Omega^{k+l} CS_{\rm c.c}(\R^n)$; let  $$\Omega CS_{\rm c.c}(\R^n):=
\bigoplus_{k=0}^\infty  \Omega^k CS_{\rm c.c}(\R^n)$$ stand for the $\N_0$ graded
algebra (also filtered by  the symbol order) of classical symbol valued forms on $\R^n$ with
constant coefficients.\\
We  also consider the sets   $\Omega^k CS_{\rm c.c}^{\Z}(\R^n):=  \bigcup_{a\in\Z}
 \Omega^k\,CS^a_{\rm c.c}(\R^n) $ of integer order classical symbols valued $k$-forms  and
$\Omega^k CS^{\notin\, \Z}(U):=\bigcup_{a\notin\, \Z} \Omega^k\,CS_{\rm
  c.c}^a(U) $ of non integer order classical symbol valued $k$-forms. Clearly,
$\Omega CS_{\rm c.c}^\Z (\R^n):=
\bigoplus_{k=0}^\infty  \Omega^k CS_{\rm c.c}^\Z(\R^n)$ is  a subalgebra
of $\Omega CS_{\rm c.c}^\Z (\R^n)$.
\\
\begin{defn}
Let  ${\cal S}\subset CS_{\rm c.c}(\R^n)$ be a set containing smoothing symbols.   We call a  linear form\footnote{By linear we mean that $\rho(\alpha_1\,
    \sigma_1+\alpha_2\, \sigma_2)= \alpha_1\, \rho(\sigma_1)+\alpha_2\,
    \rho(\sigma_2)$ whenever $\sigma_1, \sigma_2, \alpha_1\,
    \sigma_1+\alpha_2\, \sigma_2 $ lie in ${\cal S}$.}   $\rho: {\cal S}
  \to \C$ {\rm singular} if it vanishes on smoothing symbols, and {\rm
    regular} otherwise.
\end{defn}
A linear form $\rho: {\cal S}\to \C$  extends to a linear form  $\tilde \rho: \Omega{\cal S}\to \C$ defined by
$$\tilde\rho\left(\alpha(\xi)\, d\xi_{i_1}\wedge\cdots
\wedge  d\xi_{i_k}\right):=\rho(\alpha)\,  \delta_{k-n},$$ with
$i_1<\cdots<i_k$. Here we have set
$$\Omega^k{\cal S}:=\{\sum_{\vert I\vert = k} \alpha_I (\xi)\, d\xi_I,  \quad
\alpha_I\in {\cal S}\}.$$   Exterior differentiation on forms
extends to  symbol valued forms (see (5.14) in \cite{LP}):
\begin{eqnarray*}
d: \Omega^k CS_{\rm c.c}(\R^n)&\to & \Omega^{k+1} CS_{\rm c.c}(\R^n)\\
\alpha(\xi)\, d\xi_{i_1}\wedge\cdots \wedge d\xi_{i_k}&\mapsto &\sum_{i=1}^n
\partial_i\alpha(\xi)\, d\xi_i\wedge d\xi_{i_1}\wedge\cdots\wedge d\xi_{i_k}.
\end{eqnarray*}
We call a  symbol valued form $\alpha$ closed if  $d\alpha=0$ and exact
if  $\alpha=d\,\beta$ where $\beta$ is a symbol valued form; this
gives rise to the following  cohomology groups
$$H^kCS_{\rm c.c}(\R^n):=\{\alpha \in \Omega^k CS_{\rm c.c}(\R^n), \quad
d\alpha=0\}\, /\, \{d\, \beta, \beta \in \Omega^{k-1} CS_{\rm c.c}(\R^n)\}.$$
We call a  symbol valued form $\alpha$ closed ``up to a smoothing symbol''  if
$d\alpha\sim 0$ and exact
``up to a smoothing symbol'' if   $\alpha\sim d\,\beta$ where $\beta$ is a
symbol valued form. Since $\alpha\sim d \beta \Rightarrow d\alpha\sim 0$, this
gives rise to the following  cohomology groups
$$H_\sim^kCS_{\rm c.c}(\R^n):=\{\alpha \in \Omega^k CS_{\rm c.c}(\R^n), \quad
d\alpha\sim 0\}\, /\, \{\alpha\sim d\, \beta, \beta \in \Omega^{k-1} CS_{\rm
  c.c}(\R^n)\}.$$
The next two paragraphs  are dedicated to the description of the set of
top degree
forms which are exact ``up to smoothing operators'' (see Corollary
\ref{cor:resexact}). The uniqueness of the residue as a closed singular linear
form on the algebra of symbols then follows  (see Theorem
\ref{thm:uniquenessres}).

\subsection{Stokes' property versus translation invariance}
\begin{lem}\label{lem:closednessconditions}Let  $\rho:{\cal S}\subset CS_{\rm c.c}(\R^n)\to \C$ be a linear
  form. The following two conditions are equivalent:
\begin{eqnarray*}
\exists j\in \{1, \cdots, n\}\quad {\rm s.t.}\quad  \sigma=\partial_{j}
\tau \in {\cal S}&\Longrightarrow& \rho (\sigma)=0\\
\alpha=  d\, \beta \,   \in \Omega^n{\cal S}&\Longrightarrow &\widetilde \rho
(\alpha)=0
\end{eqnarray*}
\end{lem}
{\bf Proof:} We first show that the second condition follows from the first one. Since $\tilde \rho$ vanishes on forms of degree $<n$ we can
assume that  $\alpha$ is a homogeneous form of degree $n$ and show that the first condition implies that $\tilde\rho(\alpha)=0$.
Write $\alpha= d\left( \sum_{\vert J\vert = n-1} \beta_J \,
  d\xi_{i_1}\wedge \cdots\wedge d\xi_{i_{J}}\right)= \sum_{i=1}^n  \sum_{\vert
  J\vert= n-1} \partial_i  \beta_J \,
  d\xi_i\wedge d\xi_{i_1}\wedge \cdots\wedge d\xi_{i_{n-1}}$ then
$\tilde \rho( \alpha)=   \sum_{i=1}^n  \sum_{\vert
  J\vert = n-1}\rho( \partial_i  \beta_{J})  $ vanishes by the
  first condition. \\ Conversely, if $\sigma=\partial_i \tau$ then  \begin{eqnarray*} \alpha&=&\sigma(\xi)\, d\xi_{1}\wedge \cdots
\wedge d\xi_{n}\\
&=& \partial_i\tau(\xi)\, d\xi_{1}\wedge \cdots
\wedge d\xi_{n}\\
&=&  (-1)^{i-1}\, d\left(\tau(\xi)\,
  d\xi_{1}\wedge \cdots \wedge d\xi_{i-1} \wedge d\hat \xi_i \wedge d\xi_{i+1} \wedge  \cdots
\wedge d\xi_{n}\right)\\
&=& d\left(   (-1)^{i-1}\, \tau(\xi)\,\
  d\xi_{1}\wedge \cdots \wedge d\xi_{i-1} \wedge d\hat \xi_i \wedge d\xi_{i+1} \wedge  \cdots
\wedge d\xi_{n}\right)
\end{eqnarray*} is an exact form $\alpha=d\, \beta_i$ where we have set
$\beta_i:=  (-1)^{i-1}\,\tau_i(\xi)\,
  d\xi_{1}\wedge \cdots \wedge d\xi_{i-1} \wedge d\hat \xi_i \wedge d\xi_{i+1} \wedge  \cdots
\wedge d\xi_{n}$. If the  second condition is satisfied then
 $\widetilde \rho\circ d(\beta_i)=0$ from which the first condition
$ \rho\circ\partial_i(\tau)=0$ follows. \endsquare\\ \\
Following the terminology used in noncomutative geometry, we set the following definitions.
\begin{defn}
 A  linear form $\tilde\rho:\Omega {\cal S}\subset \Omega CS_{\rm
   c.c}(\R^n)\to \C$  is {\it  closed
} when it satisfies  the
equivalent conditions of Lemma \ref{lem:closednessconditions}. We also say by
extension that $\rho$ is closed if $\tilde\rho$ is or with
the analogy with the ordinary integral in mind, that
$\rho$ satisfies {\it Stokes' property}.
\end{defn}
\begin{rk}\begin{enumerate}
\item A closed linear form $\tilde \rho: \Omega CS_{\rm c.c}(\R^n)\to \C$ clearly induces
a linear form $H^\bullet CS_{\rm c.c}(\R^n)\to \C$.
\item  When $\rho$ is singular, closedness of $\tilde \rho$ is equivalent
  to the fact that $$\alpha\sim d\beta\Longrightarrow \tilde \rho(\alpha)=0.$$
 A closed singular linear form  therefore induces a  linear form 
 $ H_\sim^\bullet CS_{\rm c.c}(\R^n)\to \C$.
\end{enumerate}
\end{rk}
 Closedness turns out to be equivalent   to translation
invariance for any linear map on classical symbols which fulfills Stokes' property
on symbols of negative enough order.   This extends results  of \cite{MMP} .
\begin{prop}\label{prop:Stokestransl} Let ${\cal S}\subset  CS_{\rm c.c}(\R^n)$ be a set stable under
  translations and derivatives such that there is some positive  integer $K$
 $$CS_{\rm c.c}^{<-K}(\R^n) \subset {\cal S}.$$
 Let  $\rho:{\cal S} \to \C$ be a linear map with the Stokes' property on
 $CS_{\rm c.c}^{<-K}(\R^n)$ i.e.
$$\exists j\in \{1, \cdots, n\}\quad {\rm s.t.}\quad  \sigma=\partial_{j}
\tau \in{\cal S}\cap CS_{\rm c.c}^{<-K}(\R^n)\Longrightarrow \rho(\sigma)=0.$$
Then  for any $\sigma\in {\cal S}$
we have
\begin{eqnarray*}
\rho\left(\partial_{j}\sigma\right)&=&0 \quad \forall j\in \{1, \cdots, n\} \quad
{\rm (closedness}\quad {\rm  condition)}\\
\Longleftrightarrow\quad \rho\left( t_\eta^*\sigma\right)&=&\rho(\sigma)\quad \forall \eta \in \R^n \quad {\rm (translation}\quad {\rm invariance)}.
\end{eqnarray*}
\end{prop}
{\bf Proof:}  The proof  borrows  ideas from \cite{MMP}.\\
Let $\sigma\in CS_{c.c}(\R^n)$ and let us write a Taylor expansion of the map $t_\eta^*\sigma$ in a neighborhood of $\eta=0$. There is some $\theta\in ]0, 1[$ such that
$$t_\eta^*\sigma= \sum_{\vert \alpha \vert =0}^{N-1} \partial^\alpha \sigma \,\frac{ \eta^\alpha }{\alpha!}+ \sum_{\vert \alpha \vert =N} \partial^\alpha t^*_{\theta \eta}(\sigma) \,\frac{ \eta^\alpha }{\alpha!}.$$
Choosing $N$ large enough for $ \partial^\alpha t_{\theta
   \eta}(\sigma)$ to be of order $<-n$, it  follows from the  linearity of $\rho$   that
\begin{eqnarray}\label{eq:dtransl}
\rho(t_\eta^*\sigma)&=& \sum_{\vert \alpha \vert =0}^{N-1}\rho( \partial^\alpha \sigma )\,\frac{ \eta^\alpha }{\alpha!}+ \sum_{\vert \alpha \vert =N}\rho( \partial^\alpha t^*_{\theta \eta}(\sigma)) \,\frac{ \eta^\alpha }{\alpha!}\nonumber\\
&=& \sum_{\vert \alpha \vert =0}^{N-1}\rho( \partial^\alpha \sigma )\,\frac{
  \eta^\alpha }{\alpha!}+ \sum_{\vert \alpha \vert =N}\rho( \partial_{j} \,
\partial^\beta t^*_{\theta \eta}(\sigma)) \,\frac{ \eta^\alpha
}{\alpha!}\nonumber\\
&=& \sum_{\vert \alpha \vert =0}^{N-1}\rho( \partial^\alpha \sigma )\,\frac{
  \eta^\alpha }{\alpha!}
\end{eqnarray}
so that
$$\rho(t_\eta^*\sigma)-\rho(\sigma)= \sum_{\vert \alpha \vert =1}^{N-1}\rho( \partial^\alpha \sigma )\,\frac{
  \eta^\alpha }{\alpha!}$$ from which the result follows.\\
Here we set $\partial^\alpha=\partial_{j}\circ \partial^\beta$
for some multiindex  $\beta$ whenever  $\vert \alpha\vert \neq 0$ and,
choosing $N$ large enough so that the remainder term is of order $<-K$ we used the assumption that  $\rho$ verifies Stokes' property on $CS_{\rm c.c}^{<-K}(\R^n)$.
\endsquare
\subsection{A characterisation of the noncommutative residue and its kernel}
We show that the noncommutative residue is the unique singular linear form on
classical symbols on $\R^n$ with constant coefficients which fulfills  Stokes'
property. This is based on results of \cite{FGLS} and \cite{G1} (see
also \cite{L} for a generalisation to logarithmic powers) as well as results
of \cite{MMP}.\\ \\ {\it We henceforth and throughout the paper assume that  the dimension $n$ is larger
  or equal  two.}\\
\begin{defn}  The noncommutative residue  is a linear form on  $CS_{\rm c.c}(\R^n)$ defined by
$${\rm res}(\sigma):=\frac{1}{\sqrt{2\pi}^n}\,  \int_{S^{n-1}}\sigma_{-n}(\xi)\,
d\mu_S(\xi)$$
where  $$d \mu_S( \xi):= \sum_{j=1}^n (-1)^{j-1} \,
\xi_j\,d\xi_1\wedge \cdots \wedge d \hat \xi_j\wedge\cdots \wedge d\xi_n$$
denotes     the  volume measure on $S^{n-1}$ induced by the
canonical measure on $\R^n$.
\end{defn}
The noncommutative residue fulfills  Stokes'  property.
\begin{prop}\label{prop:closedres} \cite{MMP} (see also \cite{LP}) The
  noncommutative residue  vanishes on symbols  which are 
   partial
  derivatives in
  $CS_{\rm c.c}(\R^n)$ up to some smoothing
  operator:
$$\sigma\sim\partial_i \tau\Longrightarrow {\rm res}\circ \sigma= 0\quad \forall i=1, \cdots, n, \quad \forall
\sigma\in CS_{\rm c.c}(\R^n).$$
Equivalently, its extension $\widetilde {\rm res}$ to classical symbol valued
forms on $\R^n$ is closed.
\end{prop}
{\bf Proof:} Assume that $\sigma\sim \partial_i\tau$.  Since res vanishes on smoothing symbols, we can assume that
  $\sigma =\partial_i \tau$ for some $\tau\in CS_{\rm c.c}(\R^n)$ then $\sigma_{-n}=\partial_i \tau_{-n+1}$. \\
We have   $d\mu_S(\xi)= \iota_X(\Omega)(\xi)$
where $\Omega(\xi):= d\xi_1\wedge \cdots\wedge d\xi_n$ is the volume form on $\R^n$
and $X:= \sum_{i=1}^n \xi_i \partial_i$ is the Liouville field on
$\R^n$. Since the map $\xi\mapsto \sigma_{-n+1}(\xi)\, d\xi_1\wedge \cdots\wedge  \hat d\xi_i
\wedge \cdots\wedge  d\xi_n$ (where $\hat d\xi_i$ means we have omitted the
variable $\xi_i$)  is invariant under  $\xi\mapsto
t\, \xi$ for any $t>0$, we have ${\cal L}_X(  \sigma_{-n+1}(\xi)\, d\xi_1\wedge \cdots\wedge  \hat d\xi_i
\wedge \cdots\wedge  d\xi_n)=0$. Using Cartan's formula ${\cal L}_X= d\circ \iota_X+\iota_X\circ d$ we write
\begin{eqnarray*}
{\rm res}(\sigma)&=&\int_{S^{n-1}}\sigma_{-n}(\xi)\,  \iota_X(\Omega)(\xi)\\
&=&\int_{S^{n-1}}  \iota_X(\partial_i \tau_{-n+1}(\xi)\,\Omega)(\xi)\\
&=&(-1)^{i-1}\, \int_{S^{n-1}}  \iota_X\circ d ( \tau_{-n+1}(\xi)\,d\xi_1\wedge \cdots\wedge  \hat d\xi_i
\wedge \cdots\wedge  d\xi_n)(\xi)\\
&=&(-1)^{i}\, \int_{S^{n-1}} d\circ  \iota_X ( \tau_{-n+1}(\xi)\,d\xi_1\wedge \cdots\wedge  \hat d\xi_i
\wedge \cdots\wedge  d\xi_n)(\xi)\\
&=&0,
\end{eqnarray*}
where we have used the ordinary Stokes' formula in the last
equality. \endsquare\\ \\
The description of homogeneous symbols as  sum of partial derivatives induces
a similar description ``up to smoothing symbols'' for  all classical symbols with vanishing
residue. The following elementary result is very useful for that purpose.
\begin{lem}\label{lem:Euler}(Euler's
  theorem) For any homogeneous functions $f$ of degree $a$ on  $\R^n-\{0\}$
$$\sum_{i=1}^n \xi_i \partial_i f= a \, f.$$
\end{lem}
{\bf Proof:}
$$ \sum_{i=1}^n \partial_i(f( \xi))\, \xi_i={\frac{\partial}{\partial_t}}_{\vert_{t=1}}\,  f(t\, \xi)=
{\frac{\partial}{\partial_t}}_{\vert_{t=1}}\, t^a f(\xi)= a\, f(\xi).$$
\endsquare\\  The following proposition collects results from \cite{FGLS} (see
Lemma 1.3).
 \begin{prop}\label{prop:Kerres}
Any  symbol $\sigma\in CS^a_{\rm c.c}(\R^n)$ with vanishing residue
$${\rm res}(\sigma)=\int_{S^{n-1}}\sigma_{-n}(\xi)\, d\xi=0$$ is up to some
smoothing symbol,
  a finite sum  of partial
  derivatives, i.e. there exist symbols  $\tau_i\in CS^{a+1}_{\rm c.c}(\R^n), i=1,
  \cdots, n$  such that
 \begin{equation}\label{eqsigmapartial}\sigma \sim \sum_{i=1}^n
   \partial_i\tau_i. \end{equation}
In particular, given a  linear form $\rho: CS_{\rm c.c}(\R^n)\to \C$,
\begin{equation}\label{eq:rhoStokes} \rho\quad {\rm is}\quad{\rm
    singular}\quad{\rm and} \quad {\rm satisfies }\quad {\rm
    Stokes'}\quad{\rm property} \Longrightarrow {\rm
    Ker}( {\rm
  res})\subset {\rm Ker}( \rho).
\end{equation}
\end{prop}
{\bf Proof:} Equation (\ref{eq:rhoStokes}) clearly follows from equation
 (\ref{eqsigmapartial}) since $\rho$ is assumed to vanish on smoothing
 symbols. \\ To prove  (\ref{eqsigmapartial})  we write $\sigma\sim \sum_{j=0}^\infty\chi\,  \sigma_{a-j}$ with
$\sigma_{a-j}\in \Ci(\R^n-\{0\})$ homogeneous of degree $a-j$.
\begin{itemize}
\item  If $a-j\neq -n$ it  follows from Lemma
\ref{lem:Euler} that the homogeneous
  function  $\tau_{i,a-j+1}=\frac{\xi_i\, \sigma_{a-j}}{a+n-j}$ is such that $\sum_{i=1}^n
  \partial_i \tau_{i,a-j+1}=\sigma_{a-j}$ since
$ \sum_{i=1}^n \partial_i(\sigma_{a-j}(\xi)\, \xi_i)= (a+n-j)\, \sigma_{a-j}(\xi).$
\item We now consider the case  $a-j=-n$.  In polar coordinates $(r,
  \omega)\in \R^+_0\times S^{n-1}$ the Laplacian reads
$\Delta=-\sum_{i=1}^n\partial_i^2= -r^{1-n} \partial_r (r^{n-1} \partial_r)+
r^{-2}\Delta_{S^{n-1}}$. Since
$ \Delta( f(\omega) r^{2-n})=r^{-n}\, \Delta_{S^{n-1}} $ we have
$  \Delta( f(\omega) r^{2-n})= \sigma_{-n}(r\omega)\Leftrightarrow
\Delta_{S^{n-1}}f=\left(\sigma_{-n}\right)_{S^{n-1}}.$  Setting $F(r\,
\omega):= f(\omega)\, r^{2-n}$ it follows that the equation
$\Delta F=\sigma_{-n} $  has a solution  if and only if $\sigma_{-n}\in {\rm Ker} \Delta_{S^{n-1}}^{\perp_{S^{n-1}}}$ i.e. if
res$(\sigma)=0$. In that case, $\sigma_{-n}= \sum_{i=1}^n \partial_i
\tau_{i,-n+1}$ where we have set  $\tau_{i,-n+1}:= \partial_iF$.
\end{itemize}
Let $\tau_i\sim \sum_{j=1}^\infty \chi\, \tau_{i,a-j+1}$ then
\begin{equation}\label{eq:sigmapartialtau} \sigma\sim \sum_{i=0}^n\sum_{j=0}^\infty \chi\, \partial_i \tau_{i,a-j+1}\sim  \sum_{i=1}^n
\partial_i \tau_i\end{equation}since $\partial_i\chi$ has compact support so that the difference
$ \sigma- \sum_{i=1}^n \partial_i\tau_i$ is smoothing. Since the $\tau_i$ are
by construction of order $a+1$,  statement (\ref{eqsigmapartial}) of the proposition
follows.
\endsquare
\\ \\ The following proposition gives a characterisation of the kernel of the
noncommutative residue.
\begin{cor} \label{cor:resexact}Top degree  symbol valued forms which are
  exact up to smoothing symbols coincide with forms with vanishing
(extended)  residue:
\begin{equation}\label{eq:exactres}{\rm
  Ker}\left(\widetilde {\rm res}_{\vert_{\Omega^{n} CS_{\rm
        c.c}(\R^n)}}\right)=\{\alpha\sim d\beta, \quad\beta\in
\Omega^{n-1} CS_{\rm c.c}(\R^n)\}
\end{equation} which implies that
$$H_\sim^nCS_{\rm c.c}(\R^n):=\{\alpha \in \Omega^n CS_{\rm c.c}(\R^n), \quad
d\alpha\sim 0\}\, /\,{\rm
  Ker}\left(\widetilde {\rm res}_{\vert_{\Omega^{n} CS_{\rm
        c.c}(\R^n)}}\right) .$$
Given  a linear
form $\rho: CS_{\rm c.c}(\R^n)\to \C$
\begin{equation}\label{eq:rhoclosed}\tilde \rho \quad {\rm is}\quad {\rm
    closed }\quad {\rm and}\quad{\rm
    singular} \Longleftrightarrow {\rm
    Ker}(\widetilde {\rm
  res}) \subset {\rm Ker}(\tilde \rho).
\end{equation}
\end{cor}
{\bf Proof:} Equation  (\ref{eq:rhoclosed}) clearly follows from the first part of the
assertion.\\
Let us turn to the first part of the assertion and prove
(\ref{eq:exactres}).
By Proposition \ref{prop:closedres}, we  know that exact forms  lie in the
kernel of the residue up to smoothing
symbols. \\ To prove the other inclusion, let
$$\alpha= \sum_{\vert J\vert=0}^k \alpha_J(\xi)\, d\xi_{i_1}\wedge \cdots
\wedge d\xi_{i_{\vert J\vert}}$$
(we can choose $i_1<\cdots <i_{\vert J\vert}$ without loss of generality)  has  vanishing residue $\widetilde {\rm
  res}(\alpha)=0$. Then either $\vert J\vert <n$  or $\vert J\vert=n$ in which
case $i_1=1, \cdots, i_n=n$ and  $\widetilde{\rm res}(\alpha)={\rm res}(\alpha_n)=0$. In
that case, we can apply Proposition \ref{prop:Kerres} to $\sigma:=\alpha_n$
and write:
\begin{eqnarray*} \sigma(\xi)\, d\xi_{1}\wedge \cdots
\wedge d\xi_{n}
&\sim& \sum_{i=1}^n \partial_i\tau_i(\xi)\, d\xi_{1}\wedge \cdots
\wedge d\xi_{n}\\
&\sim& \sum_{i=1}^n  (-1)^{i-1}\, d\left(\tau_i(\xi)\,\wedge
  d\xi_{1}\wedge \cdots \wedge d\xi_{i-1} \wedge d\hat \xi_i \wedge d\xi_{i+1} \wedge  \cdots
\wedge d\xi_{n}\right)\\
&\sim& d\left(  \sum_{i=1}^n  (-1)^{i-1}\, \tau_i(\xi)\,\wedge
  d\xi_{1}\wedge \cdots \wedge d\xi_{i-1} \wedge d\hat \xi_i \wedge d\xi_{i+1} \wedge  \cdots
\wedge d\xi_{n}\right)
\end{eqnarray*} which proves (\ref{eq:exactres}).\endsquare
\begin{thm}\label{thm:uniquenessres}
Any singular linear form  $\rho: CS_{\rm c.c}(\R^n)\to \C$ with Stokes' property
 is proportional to the residue, i.e. $\rho=c\cdot {\rm res}$ for some
 constant $c$.\\ Equivalently, any closed singular linear form  $\widetilde \rho:
 \Omega CS_{\rm c.c}(\R^n)\to \C$  is proportional to the residue extended to
 forms, i.e. $\widetilde \rho=c\cdot \widetilde {\rm res}$ for some
 constant $c$.
\end{thm}
{\bf Proof:} By Proposition \ref{prop:closedres},   $ \rho$ satisfies Stokes' property
implies that  $\rho$ vanishes on Ker (res). \\
Since  $\sigma-{\rm res} (\sigma) \,\frac{\sqrt{2\pi}^n}{{\rm Vol}(S^{n-1})} \,
\vert \xi\vert^{-n}\, \chi(\xi)$  has vanishing residue\footnote{Here as before $\chi$ is a smooth cut-off
function which vanishes in a neighborhood of $0$ and is identically $1$
outside the unit ball.}, we infer
that $ \rho(\sigma)=  {\rm res} (\sigma) \,\frac{\sqrt{2\pi}^n}{{\rm Vol}(S^{n-1})}
\, \rho(\xi\mapsto \vert \xi\vert^{-n}\, \chi(\xi))$ from which the  statement of the theorem
follows setting $c:= \frac{\sqrt{2\pi}^n}{{\rm Vol}(S^{n-1})}
\, \rho(\xi\mapsto \vert \xi\vert^{-n}\, \chi(\xi))$. Since $\rho$ vanishes on
smoothing symbols by assumption, this constant is independent of the choice of $\chi$.
\endsquare
\subsection{ A characterisation of the cut-off regularised integral
  $\cutoffint_{\R^n}$ in terms of Stokes' property}
Let us recall the construction of a useful extension of the ordinary integral
given by the cut-off regularised integral.\\
  For any $R>0$,  $B(0, R)$
denotes the ball of radius $R$ centered at $0$ in $\R^n$.
We recall that given a symbol $\sigma\in CS^a_{\rm c.c}(\R^n)$,  the map
$R\mapsto \int_{B(0, R)} \sigma (\xi)\,d\,  \xi$ has an asymptotic
expansion as $R\to\infty$ of  the form (here we use the notations of
(\ref{eq:asymptsymb}):$$\int_{B(0, R)}\sigma(\xi) \, d\, \xi\\
\sim_{ R\to \infty}\alpha_0(\sigma)+
\sum_{j=0,a-j+n\neq 0}^\infty  \sigma_{a-j }\,  R^{a-j+n}+
{\rm res}(\sigma) \cdot  \log  R$$  where
$\alpha_0(\sigma)$ is the constant term  given by:
\begin{eqnarray}\label{eq:constanttermclassical}
\cutoffint_{\R^n}  \sigma (\xi) \, d\,\xi& :=&  \int_{\R^n}  \sigma_{(N)} (\xi)\, \dbar\xi+\sum_{j=0}^{N-1}  \int_{B(0, 1)}
 \chi(\xi)\,   \sigma_{a-j}(\xi) \, d\, \xi\nonumber\\
&  -& \sum_{j=0, a-j+n\neq 0}^{N-1}  \frac{ 1}{a-j+n}
 \int_{S^{n-1}} \sigma_{a-j } (\omega)\,d  \mu_S(\omega).
\end{eqnarray}
This cut-off integral $\cutoffint_{\R^n}$  defines a linear form on
$CS_{\rm
  c.c}(\R^n)$ which extends the ordinary integral in the following sense; if
$\sigma$ has complex order with real part smaller than $-n$ then $\int_{B(0, R)} \sigma(\xi)\,
d\, \xi$ converges as $R\to \infty$ and $$\cutoffint_{\R^n} \sigma(\xi)\,
d\,\xi=\int_{ \R^n}  \sigma (\xi)\, d\,\xi.$$
As it is the custom for the ordinary integral we use the same symbol
$\cutoffint_{\R^n}$ for its extension to forms so that:
$$\cutoffint_{\R^n} \alpha(\xi)\, d\xi_{i_1}\wedge\cdots
\wedge  d\xi_{i_k}:=\left(\cutoffint_{\R^n}\alpha\right)\,  \delta_{k-n},$$
where we have assumed that $i_1<\cdots <i_k$
\begin{rk} Since the cut-off regularised integral $\cutoffint_{\R^n}$
  coincides on symbols of order $<-n$ with the ordinary integral which
  vanishes on partial derivatives, $\rho:=\cutoffint_{\R^n}$ fulfills
  the assumptions of Proposition  \ref{prop:Stokestransl} with ${\cal S}=
  CS_{\rm c.c}(\R^n)$. Consequently,  translation invariance of
  $\cutoffint_{\R^n}$ is equivalent to closedness:
\begin{eqnarray*}
\cutoffint_{\R^n} \partial_{j}\sigma&=&0 \quad \forall j\in \{1, \cdots, n\} \quad
{\rm (closedness}\quad {\rm  condition)}\\
\Longleftrightarrow \cutoffint_{\R^n}  t_\eta^*(\sigma)&=&\cutoffint_{\R^n}  \sigma\quad \forall \eta \in \R^n \quad {\rm
  (translation}\quad {\rm invariance)}.
\end{eqnarray*}
\end{rk}
We investigate its closedness: unfortunately,
the cut-off regularised integral does not in general satisfy Stokes' property.
\begin{prop}\cite{MMP}\label{prop:Stokesint}
For  any  $\tau\in  CS_{\rm c.c}(\R^n) $ we have
$$\cutoffint_{ \R^n} \partial_i \tau( \xi) \, d\xi=
(-1)^{i-1}\int_{ \vert\xi\vert=1} \tau_{-n+1}( \xi) \,  d\xi_1\wedge \cdots \wedge
\hat{d\xi_i}\wedge \cdots \wedge d\xi_n.$$
\end{prop}
{\bf Proof:}
\begin{eqnarray*}
\cutoffint_{\R^n} \partial_i\tau( \xi) \,  d\xi&=& {\rm fp }_{R\to \infty}
\int_{B(0, R)} \partial_i\tau( \xi) \,  d\xi\\
&=& {\rm fp }_{R\to \infty}R^n
\int_{B(0, 1)} \partial_i\tau(R \xi) \,  d\xi\\
&=& {\rm fp }_{R\to \infty}R^{n-1}
\int_{B(0, 1)} \partial_i(\tau(R \xi)) \,  d\xi\\
&=& (-1)^{i-1} {\rm fp }_{R\to \infty}R^{n-1}
\int_{B(0, 1)} d\left(\tau(R \xi) \,  \,  d\xi_1\wedge \cdots \wedge
  \hat{d\xi_i}\wedge \cdots \wedge d\xi_n\right) \\
&=& (-1)^{i-1} {\rm fp }_{R\to \infty}R^{n-1}
\int_{ \vert\xi\vert=1} \tau(R \xi) \,  \,  d\xi_1\wedge \cdots \wedge
  \hat{d\xi_i}\wedge \cdots \wedge d\xi_n \\
&=&(-1)^{i-1}  \int_{ \vert\xi\vert=1} \tau_{-n+1}( \xi) \,
d\xi_1\wedge \cdots \wedge \hat{d\xi_i}\wedge \cdots \wedge d\xi_n
\end{eqnarray*}
in view of (\ref{eq:constanttermclassical}). \endsquare\\ \\
However, the cut-off regularised integral does obey Stokes' property on specific types
of
symbols.
\begin{cor}\label{cor:examples}We have
$$\sigma(\xi)= \partial_i \tau( \xi)\Rightarrow \cutoffint_{ \R^n}  \sigma( \xi) \, d\xi=0\quad \forall i\in \{1,
\cdots, n\}$$ in the following cases:
\begin{enumerate}
\item if $\sigma$ has non integer order,
 \item if $\sigma$ has integer order $a$ and $\sigma_{a-j}(-\xi)=
      (-1)^{a-j}\sigma(\xi)\quad \forall j\in \N_0$
in odd dimension
\item  if  has integer order $a$  and  $\sigma_{a-j}(-\xi)=
      (-1)^{a-j+1}\sigma(\xi)\quad \forall j\in \N_0$   in even dimension.
\end{enumerate}

\end{cor}
{\bf Proof:} By Proposition
\ref{prop:Stokesint}, $$\cutoffint_{  \R^n} \sigma( \xi) \, d\xi=
(-1)^{i-1}\int_{ \vert \xi\vert=1} \tau_{1-n}( \xi) \,  d\xi_1\wedge \cdots \wedge \hat{d\xi_i}\wedge \cdots \wedge d\xi_n.$$
\begin{enumerate} \item If $\sigma$ has non integer order, then  so has $\tau$  which implies that  $\tau_{1-n}=0$
  so that $\cutoffint_{ R^n} \sigma( \xi) \, d\xi$ vanishes.
\item  For any holomorphic family \footnote{We refer the reader to Section 2  for the
    notion of holomorphic family of symbols.} $\sigma(z)$ in $CS_{\rm
    c.c}(\R^n)$ with non constant affine holomorphic order
$\alpha(z)$  and  such that
$\sigma(0)=\sigma$  we have by \cite{PS} (see (\ref{eq:PSsymb}) in Section 2)
$${\rm fp}_{z=0}\cutoffint_{\R^n} \sigma(z)=\cutoffint_{\R^n}
\sigma-\frac{1}{\alpha^\prime(0)} \int_{S^{n-1}}\left(
  \partial_z\sigma_{\vert_{z=0}}\right)_ {-n}\, d\mu_S.$$
We apply this to $\sigma(z)= \partial_i(\tau(z))$ with  $\tau(z)(x)=\chi(x)\,
\tau(x)\, \vert \xi\vert^{-z}$ for some smooth cut-off funciton $\chi$ which
vanishes in a neighborhood of $0$ and is identically one outside the open unit ball.
   By the first part of the corollary, since $\sigma(z)$ has
   non integer order outside a discrete set of complex numbers (which
   correspond to the poles of $\cutoffint_{\R^n}\sigma(z)$)  we have $\cutoffint_{\R^n}
\sigma(z)=\cutoffint_{\R^n} \partial_i(\tau(z))=0$ as a meromorphic map. On the
other hand, since $\left(\partial_z\sigma_{\vert_{z=0}}\right)_{\vert_{S^{n-1}}}= -\left(\tau\,\partial_i  \log \vert
    \xi\vert\right)_{\vert_{S^{n-1}}}= -\left(\tau
      \,\xi_i\right)_{\vert_{S^{n-1}}} $ it follows that $\int_{S^{n-1}}\left(
        \partial_z\sigma_{\vert_{z=0}}\right)_{-n}\, d\mu_S=\int_{S^{n-1}}\tau_{-n-1}
      \,\xi_i \, d\mu_S$. But this last quantity vanishes whenever
      $\tau_{-n-1}$ is an even function  i.e whenever $\sigma_{-n-2}=\partial_i\tau_{-n-1}$ is an
      odd function.  This holds in odd dimension if $\sigma_{a-j}(-\xi)=
      (-1)^{a-j}\sigma(\xi)$ or in even dimension if $\sigma_{a-j}(-\xi)=
      (-1)^{a-j+1}\sigma(\xi)$ so that  in both of these cases $\cutoffint_{\R^n}\sigma=0$.

\end{enumerate}\endsquare
\begin{thm}\label{thm:uniquenesscutoffint} Let ${\cal S}$ be a subset of $CS_{\rm
    c.c}(\R^n)$
  such that  $$CS_{\rm
    c.c}^{-\infty}(\R^n)\subset {\cal
    S}\subset {\rm Ker}({\rm res}).$$ Then by Proposition \ref{prop:Kerres}
\begin{eqnarray*}
\sigma\in {\cal S}\cap CS_{\rm
  c.c}^a(\R^n)&\Longrightarrow & \exists \tau_i \in CS^{a+1}_{\rm
  c.c}(\R^n)\quad{\rm s.t.}\quad  \sigma\sim \sum_{i=1}^n\partial_i\tau_i\quad \\
{\rm with}  \quad \sigma\sim\sum_{j=0}^\infty \chi\, \sigma_{a-j} \quad &{\rm and}&
\quad  \tau_i\sim \sum_{j=0}^\infty \chi\, \tau_{i,a-j+1}, i=1, \cdots, n
.
\end{eqnarray*}If for any $\sigma\in {\cal S}$ the
$\tau_i$ and  $\chi\, \tau_{i,a+1-j}, j\in \N_0$ can be chosen in ${\cal S}$ then
  any   linear form $\rho: {\cal
  S}\to \C$ which statisfies  Stokes' property   is entirely determined by its
  restriction to  $CS^{< -K}_{\rm c.c}(\R^n)$ for any  positive
  integer $K\leq n$.\\
 Equivalently, under the same conditions any  closed  linear form $\tilde\rho: \Omega{\cal
  S}\to \C$  is entirely determined by its
  restriction to  $\Omega CS^{< -K}_{\rm c.c}(\R^n)$ for any positive
  integer $K\leq n$.
\\ \\ In particular, if $\cutoffint_{\R^n}$  satisfies Stokes' property on
  ${\cal S}$ and $\rho$ is continuous\footnote{i.e. its restriction to symbols of constant order
    is continuous.}  on  ${\cal S}\cap CS^a_{\rm
    c.c}(\R^n)$ for any complex number $a$,  then there is a constant $c$ such that
$$\rho= c\cdot \cutoffint_{\R^n}.$$
\end{thm}
\begin{rk} In practice ${\cal S}$ can be described in terms of   some  condition  on the
  homogeneous components of the symbol and therefore
  automatically satisfies the requirements of the theorem.
\end{rk}
{\bf Proof:}
We write a symbol $\sigma\in CS^a_{\rm c.c}(\R^n)$
$$
\sigma = \sum_{j=0}^{N-1}\chi\,  \sigma_{a-j}+  \sigma_{(N)}$$
with $N$ any  integer chosen large enough so that $\sigma_{(N)}$ has order $<-n$.  Here $\chi$
is a smooth cut-off function which vanishes in a  neighborhood of
$0$ and is one outside the unit ball. As before, the $\sigma_{a-j}$ are
positively homogeneous of degree $a-j$.  \\
By linearity of $\rho$ we have:
\begin{equation}\label{eq:rhosigma}\rho(\sigma)= \sum_{j=0}^{N-1}\rho(\chi\,
  \sigma_{a-j})+ \rho(\sigma_{(N)}).
\end{equation}
Let now $\sigma \in {\cal S}$. Since by the  assumption on ${\cal S}$ the
symbol  $\sigma$ has vanishing
residue    we
can write as in the proof of Proposition \ref{prop:Kerres},
$\sigma_{a-j}=\partial_{i_j}\tau_{a+1-j}$ for some $i_j\in \{1, \cdots, n\}$
and some homogeneous  symbol $\tau_{a+1-j}$.\\By the closedness condition
$\rho\left(\partial_{i_j}(\chi\, \tau_{a+1-j})\right)=0$ so that $$\rho(\chi\,
\sigma_{a-j})=\rho(\chi\,
\partial_{i_j}\tau_{a+1-j})=-\rho((\partial_{i_j}\chi)\, \tau_{a+1-j}). $$
Summing over $j=1, \cdots, N-1$ we get:
\begin{equation}\label{eq:rhounique}\rho(\sigma)=  - \sum_{j=0}^{N-1}\rho\left(
(\partial_{i_j}\chi)\, \tau_{a+1-j}\right) + \rho\left( \sigma_{(N)}\right).
\end{equation}
 Another choice of primitive $\tilde\tau_{a+1-j}=\tau_{a+1-j}+c_{ j} $
modifies this expression by $c_{ j}\, \sum_{j=0}^{N-1}\rho\left(
\partial_{i_j}\chi\right)$ which vanishes.\\
Since $N$ can be chosen arbitrarily large, formula (\ref{eq:rhounique}) shows that $\rho$ is uniquely determined by its
expression on symbols of arbitrarily negative order. \\
Thus  $\rho$ is determined by its
  restriction to  $\bigcap_{K\geq n}CS_{\rm
    c.c}^{<-K}(\R^n)=CS^{-\infty}(\R^n)$. This restriction
   is continuous as a result of the continuity of the restriction of
  $\rho$ to any $CS^{a}_{\rm c.c}(\R^n)$.  Thus $\rho$ restricted to
  $CS_{\rm c.c}^{-\infty}(\R^n)$ can be seen as a tempered distribution with vanishing derivatives
  at all orders. Such a distribution  is a priori of the form
  $f\mapsto\int_{\R^n}f(\xi)\, \phi(\xi)\, d\xi$ for some smooth function
  $\phi$; since all its derivatives vanish, $\phi$ is constant so that $\rho$
  restricted to smoothing symbols
  is   proportional
 to the ordinary integral $\int_{\R^n}$\footnote{I thank E. Schrohe for drawing
   my attention to this point.}. \\ From the above discussion we
 conclude that  two closed and continuous (on symbols of constant order) linear forms $\rho_1$ and $\rho_2$ on a set ${\cal
   S}$ which satisfy the
 assumptions of the theorem are proportional.\\ The  cut-off regularised integral
 is continuous on symbols of constant order. Thus,  if it has Stokes' property on the set
 ${\cal S}$, we infer from the above uniqueness result that $\rho$ is
 proportional to $\cutoffint_{\R^n}$.
\endsquare\\ \\ Here are some examples of subsets of $CS_{\rm c.c}(\R^n)$
which fulfill the
assumptions of Theorem \ref{thm:uniquenesscutoffint} and on which the cut-off
regularised integral $\cutoffint_{\R^n}$ satisfies Stokes' property in view
of  of Corollary \ref{cor:examples}.
\begin{ex} \label{ex:nonintegerorder}The set $CS^{\notin\Z}_{\rm c.c}(\R^n)$ of non integer order
  symbols.
\end{ex}
\begin{ex} In odd dimension $n$ the set $$CS^{\rm odd}_{\rm c.c}(\R^n):= \{\sigma\in CS^\Z_{\rm c.c}(\R^n),\quad
  \sigma_{a-j}(-\xi)=(-1)^{a-j}\sigma_{a-j}(\xi)\quad \forall  \xi\in\R^n\}$$ of
  odd-class symbols.
\end{ex}
\begin{ex}In even dimension $n$ the set $$CS^{\rm even}_{\rm c.c}(\R^n):= \{\sigma\in CS^\Z_{\rm c.c}(\R^n),\quad
  \sigma_{a-j}(-\xi)=(-1)^{a-j+1}\sigma_{a-j}(\xi)\quad \forall  \xi\in\R^n\}$$ of
  even-class symbols.
\end{ex}
From these examples we get the following  straightforward application of
Theorem
\ref{thm:uniquenesscutoffint}.
\begin{cor}\label{cor:applications} Any closed linear form on $CS^{\notin\Z}_{\rm c.c}(\R^n)$, resp.
  $CS^{\rm odd}_{\rm c.c}(\R^n)$ in odd dimensions, resp. $CS^{\rm even}_{\rm
    c.c}(\R^n)$ in even dimensions is
  determined by its restriction to  symbols of arbitrarily negative order. \\ If
  it is continuous on symbols of constant order, it is proportional to the
  cut-off regularised integral $\cutoffint_{\R^n}$.
\end{cor}\vfill\eject \noindent

\section{ Existence: The canonical trace on
  odd- (resp. even-) class operators in
  odd (resp. even)  dimensions} 
We  show that the  canonical trace density ${\rm
  TR}_x(A)\, dx$ defines a global density in odd (resp. even)  dimensions for
odd-(even-) class
operators $A$ which integrates over the manifold to the (extended) canonical
trace 
$${\rm TR}(A):= \frac{1}{\sqrt{2\pi}^n} \int_M dx\, {\rm TR}_x(A).$$   To do so,
on the grounds of  results of \cite{PS}, we
carry out a continuous extension along holomorphic paths $z\mapsto A(z)\in
\Cl(M,E)$ such that $A(0)\in \Cl^{\rm odd}(M, E)$ and $A^\prime(0)\in \Cl^{\rm
  odd}(M, E)$  (resp. $A(0)\in \Cl^{\rm even}(M,
E)$  and $A^\prime(0)\in \Cl^{\rm
  even}(M, E)$), and show that the
extension is independent of the holomorphic path,  thereby extending  results of \cite{KV}
and \cite{Gr}.  Along the way we define the
noncommutative residue on the algebra of classical pseudodifferential
operators as well as the canonical trace on non integer order classical
pseudodifferential operators.
\subsection{Notations}
Let $U$ be a connected open subset of $\R^n$ where as before we assume that
$n>1$. 
 \\ For any complex number $a$, let ${\cal S}_{\rm cpt}^a(U)$ denote the set of smooth
 functions on $U\times \R^n$ called symbols with compact support in $U$, such that for any multiindices
$\beta, \gamma\in \N^n$,  there is a constant $C(\beta, \gamma)$ satisfying the following requirement:
$$\vert\partial_\xi^\beta\partial_x^\gamma \sigma(x,\xi)\vert\leq C(\beta, \gamma) \vert (1+\vert
\xi\vert)^{{\rm Re}(a)-\vert \beta\vert}$$ where Re$(a)$ stands for the real
part of $a$, $\vert
\xi\vert$ for the euclidean norm of $\xi$. We single out the subset
$CS^a_{\rm cpt}(\R^n)\subset {\cal
  S}_{\rm cpt}^a(\R^n) $ of   symbols $\sigma$, called classical symbols of
order $a$  with compact support in $U$,   such that 
\begin{equation}\label{eq:asymptsymb}
\sigma(x,\xi)=
\sum_{j=0}^{N-1} \chi(\xi)\, \sigma_{a-j}(x, \xi) +\sigma_{(N)}(x, \xi)
\end{equation}
where $\sigma_{(N)}\in {\cal S}_{\rm cpt}^{a-N}(U)$ and
 where $\chi$ is a smooth cut-off function which vanishes in a
  small
   ball of $\R^n$ centered at $0$ and which is constant equal to $1$ outside
   the unit ball. Here $\sigma_{a-j}(x, \cdot), j\in \N_0$ are positively homogeneous of degree
   $a-j$.
\\
  Let  $$CS_{\rm cpt}^{-\infty}(U)=  \bigcap_{a\in
  \C}CS^a_{\rm cpt}(U)$$ be the set of smoothing symbols with compact
support in $U$; we write $\sigma\sim \tau$ for two symbols that differ by a
smoothing symbol. \\ We equip the set $CS_{\rm cpt}^{a}(U)$ with a Fr\'echet structure with the help of the following
semi-norms labelled by multiindices $\alpha,\beta$ and integers $j\geq
0$, $N$  (see \cite{H}):
\begin{eqnarray*}
&{} & {\rm sup}_{x\in U, \xi \in \R^n} (1+\vert \xi\vert)^{-{\rm Re}(a)+\vert \beta\vert} \, \Vert \partial_x^\alpha \partial_\xi^\beta \sigma(x, \xi)\Vert;\\
&{}&  {\rm sup}_{x\in U, \xi\in \R^n}  (1+\vert \xi\vert)^{-{\rm Re}(a)+N+\vert \beta\vert}\Vert \partial_x^{\alpha} \partial_\xi^{\beta}\left(\sigma-\sum_{j=0}^{N-1} \chi(\xi)\, \sigma_{a-j}\right)(x, \xi) \Vert;\\
&{}& {\rm sup}_{x\in U, \vert\xi\vert=1}  \Vert \partial_x^{\alpha} \partial_\xi^{\beta} \sigma_{a-j}(x, \xi) \Vert.
\end{eqnarray*}
 The star product
\begin{equation}\label{eq:starproduct}\sigma\star \tau\sim\sum_{\alpha} \frac{(-i)^{\vert \alpha\vert}}{\alpha!} \partial_\xi^\alpha \sigma\,
\partial_x^\alpha \tau
\end{equation}
  of  symbols  $\sigma\in CS_{\rm cpt}^a(U)$ and $\tau\in CS_{\rm cpt}^b(U)$
  lies in  $CS_{\rm cpt}^{a+b}(U)$ provided $a-b\in \Z$.  \\
 Let  $$ CS_{\rm cpt}(U)= \langle
\bigcup_{a\in \C}CS^a_{\rm cpt}(U)\rangle $$ denote the algebra  generated by  all classical
symbols   with compact support in  $U$. We
denote by $CS_{\rm
  cpt}^{<p}(U):= \bigcup_{{\rm Re}(a)<p} CS_{\rm cpt}^a(U)$, the set of classical
symbols  of order with real part $<p$  with compact support in $U$, by  $CS_{\rm
  cpt}^{ \Z}(U):= \bigcup_{a\in\Z} CS_{\rm cpt}^a(U)$ the algebra  of  integer
order symbols, and  by   $CS_{\rm
  cpt}^{\notin \Z}(U):= \bigcup_{a\in\C- \Z} CS_{\rm cpt}^a(U)$
the set of non integer order symbols with compact support in $U$.\\
We shall also need to consider the set  introduced in \cite{KV}
 $$CS^{\rm odd}_{\rm cpt}(U):= \{\sigma\in CS^\Z_{\rm cpt}(U),\quad
  \sigma_{a-j}(-\xi)=(-1)^{a-j}\sigma_{a-j}(\xi)\quad \forall (x,\xi) \in T^* U\}$$ of
  odd-class  (also called
even-even in \cite{Gr}) symbols and the set introduced by G. Grubb (under the name even-odd)
$$CS^{\rm even}_{\rm cpt}(U):= \{\sigma\in CS^\Z_{\rm cpt}(U),\quad
  \sigma_{a-j}(-\xi)=(-1)^{a-j+1}\sigma_{a-j}(\xi)\quad \forall (x,\xi) \in
  T^* U\}$$of  even-class symbols with compact support in $U$. \\ Whereas $CS^{\rm odd}_{\rm cpt}(U)$ is
  stable under the symbol product (\ref{eq:starproduct}), $CS^{\rm even}_{\rm cpt}(U)$ is not
  since the product of two even symbols is odd. Similarly, one can check that
  the product of an odd and an even symbol is odd, two properties which 
 conflict with
  the intuition suggested by the terminology even/odd suggested by \cite{KV} 
(hence the alternative terminology used by Grubb).\\ 
 The above definitions extend to non scalar symbols. 
Given a finite dimensional vector space $V$ and any $a\in \C$  we set $$CS_{\rm cpt}^a(U,
V):= CS_{\rm cpt}^a(U)\otimes {\rm End}(V)$$  Similarly, we define $CS(U,
V)$
,  $ CS_{\rm cpt}^{ \Z}(U,V), CS_{\rm cpt}^{\notin \Z}(U, V)$ and $ CS_{\rm
  cpt}^{\rm odd}(U, V)$, $ CS_{\rm cpt}^{\rm even}(U, V)$  from 
$CS_{\rm cpt}(U), CS_{\rm cpt}^{\Z}(U), CS_{\rm cpt}^{\notin \Z}(U)$ and
$CS_{\rm cpt}^{\rm odd}(U)$,$ CS_{\rm cpt}^{\rm even}(U)$.\begin{rk}Note that $\sigma\in CS_{\rm cpt}^a(U, V)\Longrightarrow  {\rm tr}(\sigma)\in
CS^a(U)$ where tr stands for the trace on matrices.  Similar properties hold
for  $ CS_{\rm cpt}^{ \Z}(U,V), CS_{\rm cpt}^{\notin \Z}(U, V)$ and $ CS_{\rm
  cpt}^{\rm odd}(U, V)$, $ CS_{\rm cpt}^{\rm even}(U, V)$.\end{rk}
Let $M$ be an $n$-dimensional closed connected Riemannian manifold (as before $n>1$).  For $a\in
\C$, let $\Cl^{a}(M)$ denote the
linear space of classical pseudodifferential operators of order
$a$, i.e. linear maps acting on smooth functions $\Ci(M)$, which using a
partition of unity adapted to an atlas on $M$ can be written as a finite sum of operators
$$A= {\rm Op}(\sigma(A))+ R$$
where $R$ is a linear operator with smooth kernel and 
$\sigma(A)\in CS^a_{\rm cpt}(U)$ for some open subset $U\subset
\R^n$. Here we have set
$${\rm Op}(\sigma)(u):= \int_{\R^n} e^{i\langle x-y, \xi\rangle} \sigma(x,
\xi)\, u(y)\, dy\,
d\xi$$
where $\langle\cdot, \cdot\rangle$ stands for the canonical scalar product in
$\R^n$. \\The star product (\ref{eq:starproduct}) on classical symbols with
compact support induces the operator product on (properly supported) classical pseudodifferential operators since
${\rm Op}(\sigma\star \tau)={\rm Op}(\sigma) \, {\rm Op}( \tau)$. 
It follows that the product $AB$ of two classical pseudodifferential 
operators $A\in \Cl^a(M)$, $B\in \Cl^b(M)$ lies in $
\Cl^{a+b}(M)$ provided $a-b\in \Z$. Let us
denote by $\Cl(M)=\langle \bigcup_{a\in \C}\Cl^{a}(M)\rangle$ the algebra
generated by all classical pseudodifferential operators acting on $\Ci(M)$. 
\\
  Given a finite rank vector bundle $E$ over $M$ we set
$ \Cl^a(M, E):= \Cl^a(M)\otimes {\rm
  End}(E)$, $ \Cl(M, E):= \Cl(M)\otimes {\rm
  End}(E)$.
\begin{rk}  Note that  if $A\in \Cl^a(M, E)$, in a local trivialisation
$E_{\vert_U}\simeq U\times V$ over an open subset $U$ of $M$, the map  $(x,
\xi)\mapsto \sigma(A)(x, \xi)$ lies in $CS^a(U, V)$.
\end{rk}
$\Cl^a(M, E)$   inherits a
Fr\'echet structure via the Fr\'echet structure on classical symbols of order $a$.\\
The algebras $  \Cl^{\in
  \Z}(M, E)$,$ \Cl^{\notin
  \Z}(M, E)$,  $ \Cl^{\rm odd}(M, E), \Cl^{\rm even}(M, E)$ are defined
similarly using
trivialisations of $E$   from $CS^\Z_{\rm cpt}(U)$
$\Cl^a_{\rm cpt}(U), \Cl^{\notin
  \Z}_{\rm cpt}(U)$ and $ \Cl_{\rm cpt}^{\rm odd}(M)$.\\
\subsection{Classical symbol valued forms on an open subset}
The notations introduced in  paragraph 1.1 for symbols on $\R^n$ with constant
coefficients easily extend to symbols with
support in an open subset of $U\subset \R^n$ with varying coefficients.  \\
Let $U$ be a connected open subset of $\R^n$ as before.
We borrow from \cite{MMP} (see also \cite{LP}) the following notations and
some of the subsequent definitions.
\begin{defn}Let $k$ be a non negative integer,
$a$ a complex number. We let
  \begin{eqnarray*}
\Omega^k\,CS^a_{\rm cpt}(U)&= &\{\alpha \in \Omega^k(T^*U),\quad
 \alpha = \sum_{I, J\subset \{1, \cdots, n\},  \vert I
\vert + \vert J
\vert = k} \alpha_{IJ} (x,\xi) \, d\xi_I\wedge dx_J\\
 {\rm with } &{}& \quad \alpha_{IJ}  \in CS_{\rm cpt}^{a-\vert I\vert}(U)\}
\end{eqnarray*}
denote the set of order $a$ classical symbol valued forms on $U$ with
compact support.
Let  \begin{eqnarray*}
\Omega^k\,CS_{\rm cpt}(U)&= &\{\alpha \in \Omega^k(T^*U), \quad
 \alpha = \sum_{I, J\subset \{1, \cdots, n\}, \vert I
\vert  +\vert J
\vert = k} \alpha_{IJ} (x,\xi) \,  d\xi_I\wedge dx_J\\
 {\rm with } &{}& \quad \alpha_{IJ} \in CS_{\rm cpt}(U)\}
\end{eqnarray*} denote the set of classical
symbol valued $k$-forms on $U$ of all orders with compact support.
  \end{defn}
The exterior product on forms combined with the star product on symbols induces a
product  $\Omega^k CS_{\rm cpt}(U)\times \Omega^l CS_{\rm cpt}(U)\to
\Omega^{k+l} CS_{\rm cpt}(U)$; let  $$\Omega CS_{\rm cpt}(U):=
\bigoplus_{k=0}^\infty  \Omega^k CS_{\rm cpt}(U)$$ stand for the $\N_0$ graded
algebra (also filtered by the symbol order) of classical symbol valued forms on $U$ with
compact support.\\
We shall also consider the sets   $\Omega^k CS_{\rm cpt}^{\Z}(U):=  \bigcup_{a\in\Z}
 \Omega^k\,CS^a_{\rm cpt}(U) $ of integer order classical symbols valued
 $k$-forms,
$\Omega^k CS^{\notin\, \Z}(U):=\bigcup_{a\notin\, \Z} \Omega^k\,CS_{\rm
  cpt}^a(U) $ of non integer order classical symbol valued $k$-forms,
$\Omega^k CS_{\rm cpt}^{\rm odd}(U)$, resp. $\Omega^k CS_{\rm cpt}^{\rm
  even}(U)$ of odd- (resp. even-) classical symbol valued $k$-forms.
 \\ \\ Exterior differentiation on forms
extends to  symbol valued forms (see (5.14) in \cite{LP}):
\begin{eqnarray*}
d: \Omega^k CS_{\rm cpt}(U)&\to & \Omega^{k+1} CS_{\rm cpt}(U)\\
\alpha_{IJ}(x,\xi)\, d\xi_{I}\wedge  dx_{J}&\mapsto &\sum_{i=1}^{2n}
\partial_i\alpha_{IJ}(\xi)\, d u_i\wedge d\xi_{I}\wedge dx_{J},
\end{eqnarray*}
where $u_i=\xi_i, \partial_i=\partial_{\xi_i}$ if $1\leq i\leq n$ and
$u_i=x_i, \partial_i=\partial_{x_i}$ if $n+1\leq i\leq 2n$.\\
As before, we call a  symbol valued form $\alpha$ closed if  $d\alpha=0$ and exact
if  $\alpha=d\,\beta$ where $\beta$ is a symbol valued form; this
gives rise to the following  cohomology groups
$$H^kCS_{\rm cpt}(U):=\{\alpha \in \Omega^k CS_{\rm cpt}(U), \quad
d\alpha=0\}\, /\, \{d\, \beta, \quad\beta \in \Omega^{k-1} CS_{\rm
  cpt}(U)\}.$$ \\ \\
Let  ${\cal D}(U)\subset  CS_{\rm cpt}(U)$ be a set containing smoothing symbols.   We call a  linear form\footnote{By linear we mean that $\rho(\alpha_1\,
    \sigma_1+\alpha_2\, \sigma_2)= \alpha_1\, \rho(\sigma_1)+\alpha_2\,
    \rho(\sigma_2)$ whenever $\sigma_1, \sigma_2, \alpha_1\,
    \sigma_1+\alpha_2\, \sigma_2 $ lie in ${\cal D}(U)$.}   $\rho: {\cal D}(U)
  \to \C$ {\rm singular} if it vanishes on smoothing symbols, and {\rm
    regular} otherwise.
\\ A linear form $\rho: {\cal D}(U)
  \to \C$ extends to a linear form  $\tilde \rho: \Omega {\cal D}(U)\to \C$ defined by
$$\tilde\rho\left(\alpha_{IJ}(x,\xi)\, d\xi_{i_1}\wedge\cdots
\wedge  d\xi_{i_{\vert I\vert}}\wedge\, dx_{j_1}\wedge\cdots
\wedge  d  x_{j_{\vert J\vert}} \right):= \rho(\alpha_{IJ})\,  \delta_{\vert
I\vert +\vert J\vert-2n},$$ with
$i_1<\cdots<i_{\vert I\vert}$, $j_1<\cdots<j_{\vert J\vert}$. Here we have set
$$\Omega^k{\cal D}(U):=\{\sum_{\vert I\vert+\vert J\vert \leq k} \alpha_{IJ} (x,\xi)\,d\xi_I\wedge dx_J,  \quad
\alpha_{IJ}\in {\cal D}(U)\}.$$
\\ This is a straightforward generalisation of Lemma
\ref{lem:closednessconditions}.
\begin{lem}\label{lem:closednessconditionsU} Let  $\rho:{\cal D}(U)\subset CS_{\rm cpt}(U)\to \C$ be a linear
  form. The following two conditions are equivalent:
\begin{eqnarray*}
\exists i,j\in \{1, \cdots, n\}\quad {\rm s.t.}\quad  \sigma=\partial_{\xi_i}
\tau \in {\cal D}(U)\quad  {\rm or} \quad  \sigma=\partial_{x_j}
\tau \in {\cal D}(U) &\Longrightarrow& \rho (\sigma)=0\\
\alpha=  d\, \beta \,   \in \Omega^n{\cal D}(U)&\Longrightarrow &\widetilde \rho
(\alpha)=0
\end{eqnarray*}
\end{lem}
As before we call {\it closed }  a linear form $\tilde \rho$ obeying the second
condition and by extension $\rho$ is then also said to be closed. We also say
that $\rho$ satisfies Stokes' condition.
\begin{rk}A closed linear form $\tilde \rho$ on $ \Omega CS_{\rm cpt}(U)$
induces a linear form $\bar \rho: H^\bullet CS_{\rm cpt}(U)\to \C$.
\end{rk}
\begin{prop}A linear form $\rho: {\cal D}(U)\subset CS_{\rm cpt}(U)\to \C$ is
  closed whenever 
$$\rho\left(\{\sigma, \tau\}_\star\right)=0\quad \forall \sigma, \tau \in
CS_{\rm cpt}(U), \quad {\rm s.t.} \quad\{\sigma, \tau\}_\star\in  {\cal D}(U)$$
where we have set:
$$\{\sigma, \tau\}_\star:= \sum_{\alpha} \frac{(-i)^{\vert
    \alpha\vert}}{\alpha!}\left( \partial_\xi^\alpha \sigma \partial_x^\alpha \tau- 
\partial_x^\alpha \sigma \partial_\xi^\alpha \tau\right).$$
\end{prop}

{\bf Proof:}
If the linear form is closed, we can perform integration by parts and write:
\begin{eqnarray*}
\rho\left(\{\sigma, \tau\}_\star\right)&=& \sum_{\alpha} 
\frac{(-i)^{\vert\alpha\vert}}{\alpha!} \rho\left(\partial_\xi^\alpha \sigma \partial_x^\alpha \tau- 
\partial_x^\alpha \sigma \partial_\xi^\alpha \tau\right)\\
&=& \sum_{\alpha} \frac{(-i)^{\vert
    \alpha\vert}}{\alpha!} \rho\left(\partial_x^\alpha \sigma \partial_\xi^\alpha \tau- 
\partial_x^\alpha \sigma \partial_\xi^\alpha \tau\right)\\
&=&0.
\end{eqnarray*}Conversely, if the linear form vanishes on brackets $\{\cdot,
\cdot\}_\star$ contained in ${\cal D}(U)$ then for any $\sigma\in CS_{\rm cpt}(U)$ such that
$\partial_{x_i}\sigma=i \, \{ \sigma, \xi_i\}_\star\in {\cal D}(U)$ we have 
$$\rho(\partial_{x_i}\sigma)=i \,\rho\left(\{ \xi_i,\sigma\}_\star \right)= 0$$
and similarly for any 
$\sigma\in CS_{\rm cpt}(U)$ such that
$\partial_{\xi_i}\sigma= i\,\{ x_i, \sigma\}_\star\in {\cal D}(U)$ we have 
$$\rho(\partial_{\xi_i}\sigma)= i\,\rho\left(\{ x_i,\sigma\}_\star \right)= 0.$$
\endsquare

\subsection{The noncommutative residue}

\begin{defn}  The noncommutative residue of a symbol $\sigma \in CS_{\rm cpt}(U)$
  is defined by
$${\rm res}(\sigma):=\frac{1}{(2\pi)^n}\, \int_Udx \int_{S^{n-1}}\sigma_{-n}(x,\xi)\,
\mu_S(\xi)=\frac{1}{\sqrt{2\pi}^n}\,  \int_{U}{\rm res}_x(\sigma)\,d\, x  $$
where
${\rm res}_x(\sigma):=\frac{1}{\sqrt{2\pi}^n} \int_{ S^{n-1}}\sigma_{-n}(x, \xi)\, \mu_S(\xi)$
 is the residue density at point $x$ and where  as before  $$ \mu_S( \xi):= \sum_{j=1}^n (-1)^j \,
\xi_j\,d\xi_1\wedge \cdots \wedge d \hat \xi_j\wedge\cdots \wedge d\xi_n$$
denotes     the  volume measure on $S^{n-1}$ induced by the
canonical measure on $\R^n$.
\end{defn}
\begin{lem}\label{lem:resStokesU}The noncommutative residue is a
 singular closed linear form  on $CS_{\rm cpt}(U)$ which restricts  to a continuous
map on
  each $CS_{\rm cpt}^a(U)$, $a\in C$.
\end{lem}
{\bf Proof:} The continuity follows from the definition of the residue \footnote{Note that
  this continuity holds only on symbols of constant order; it breaks down
  if one lets the order vary.}. Stokes' property  follows from Stokes'
property of the ordinary   integral on $C_{\rm cpt}^\infty(U)$ combined with the fact
that the residue density ${\rm res}_x$ vanishes on derivatives
$\partial_{\xi_j}$ which follows from Proposition \ref{prop:closedres}. \endsquare
\\ \\
Using a partition of unity, one can patch up the residue on classical symbols with compact
support  to build a noncommutative residue on classical  operators on a closed
manifold $M$ introduced by Wodzicki \cite{W1} (see also \cite{G1}).
\begin{defn}The noncommutative residue of    $A\in \Cl(M,E)$ is defined
by
$${\rm res}(A):=\frac{1}{(2\pi)^n}\, \int_M dx \int_{S_x^*M}{\rm tr}_x\left(\sigma(A)\right)_{-n}(x,\xi)\,
\mu_S(\xi)=\frac{1}{\sqrt{2\pi}^n}\,  \int_{M}{\rm res}_x(A)\,d\, x  $$
where
${\rm res}_x(A):=\frac{1}{\sqrt{2\pi}^n} \int_{S_x^*M}{\rm tr}_x\left(\sigma(A)\right)_{-n}(x, \xi)\, \mu_S(\xi)$
is  the residue density at point $x$ and where  as before  $$ \mu_S( \xi):= \sum_{j=1}^n (-1)^j \,
\xi_j\,d\xi_1\wedge \cdots \wedge d \hat \xi_j\wedge\cdots \wedge d\xi_n$$
denotes     the  volume measure on the cotangent sphere $S_x^*M$ induced by the
canonical measure on the cotangent space $T_x^*M$ at point $x$. Here ${\rm
  tr}_x$ stands for the fibrewise trace on the vector bundle End$(E)$. 
\end{defn}
\begin{rk}It follows from the definition that the residue is continuous on
  each $\Cl^a(M,E)$, $a\in \C$.
\end{rk}
We derive the cyclicity of the residue on operators  from Stokes' property of the residue
on symbols.
\begin{prop}
$${\rm res}\left([A,B]\right)=0\quad \forall A, B\in \Cl(M,E) .$$
\end{prop}
{\bf Proof:}
The product of  two \pdos  $\, A, B$ in $\Cl(M,E)$
reads
\begin{equation}\label{eq:product}
A\, B=\sum_{\vert \alpha\vert \leq N}\frac{(-i)^{\vert \alpha\vert}}{\alpha!} {\rm Op}(\partial_\xi^\alpha \sigma(A)\,
\partial_x^\alpha \sigma(B))+ R_N(A\, B)
\end{equation}
  for any integer $N$  and with   $R_N(A\, B)$ of order $<a+b-N$ where $a, b$
  are the orders of $A, B$ respectively.  Hence
$$
[A, B]=\sum_{\vert \alpha\vert \leq N}\frac{(-i)^{\vert \alpha\vert}}{\alpha!} {\rm Op}\left(\partial_\xi^\alpha \sigma(A)\,
\partial_x^\alpha \sigma(B)- \partial_\xi^\alpha \sigma(B)\,
\partial_x^\alpha \sigma(A)\right)+ R_N([A, B])$$
with similar notations.\\
Applying  the noncommutative residue
 on either side, choosing $N$ such that $a+b-N<-n$ we have
\begin{eqnarray*}
{\rm res}\left([A, B]\right)&=&\sum_{\vert \alpha\vert \leq N}\frac{(-i)^{\vert
    \alpha\vert}}{\alpha!} \int_M d\,x\,
\cutoffint_{\R^n}{\rm tr}_x\left(\partial_\xi^\alpha \sigma(A)\,
\partial_x^\alpha \sigma(B)- \partial_\xi^\alpha \sigma(B)\,
\partial_x^\alpha \sigma(A)\right)\, d\,\xi+ {\rm res}\left(R_N([A,
B])\right)\\
&=& \sum_{\vert \alpha\vert \leq N}\frac{(-i)^{\vert
    \alpha\vert}}{\alpha!} \int_M d\,x\,
\cutoffint_{\R^n}{\rm tr}_x\left(\partial_\xi^\alpha \sigma(A)\,
\partial_x^\alpha \sigma(B)-
\partial_x^\alpha \sigma(A)\, \partial_\xi^\alpha \sigma(B)\right)\, d\,\xi\\
&=& 0.
\end{eqnarray*}
In the last identity we used Stokes' property for residue on symbols to
implement repeated integration by parts combined with the fact that the
residue vanishes on symbols of order $<-n$ and the cyclicity of the ordinary
trace on matrices. \endsquare
\subsection{The canonical trace on non integer order operators}
The cut-off regularised integral extends to $CS_{\rm cpt}(U)$.
\begin{defn}  For any   $\sigma \in CS_{\rm cpt}(U)$ the cut-off regularised integral  of $\sigma$
  is defined by
$$\cutoffint_{T^*U}\sigma:= \int_Udx \cutoffint_{T_x^*U}d\xi\,
\sigma(x,\xi).  $$
It extends to  pseudodifferential symbol valued forms by
$$\cutoffint_{T^*U}\alpha_{IJ} \, d\xi_I\wedge dx_J:=
\left(\cutoffint_{T^*U}\alpha_{IJ} \right)\, \delta_{\vert I\vert +\vert J\vert=2n}$$
where $d\xi_I:= d\xi_{i_1}\wedge\cdots\wedge d\xi_{i_k}$ with $i_1<\cdots <i_k$
and $dx_J:= dx_{j_1}\wedge\cdots\wedge dx_{j_l}$ with $j_1<\cdots <j_l$.
\end{defn}
\begin{lem}\label{lem:cutoffintStokesU}The cut-off regularised integral is a
  linear form  on $CS_{\rm cpt}(U)$ which restricts  to a continuous linear
form on
  each $CS_{\rm cpt}^a(U)$ and  satisfies
  Stokes' property on non integer order symbols:
$$\left( \exists j=1, \cdots, n, \quad \sigma=\partial_{x_j}\tau \quad{\rm
    or} \quad \sigma=\partial_{\xi_j}\tau\quad{\rm with }\quad \sigma \in  CS^{\notin\Z}_{\rm cpt}(U)\right) \Rightarrow
\cutoffint_{U\times\R^n}
  \sigma=0.$$
Equivalently, it extends to a linear form on $\Omega CS_{\rm cpt}(U)$ which restricts  to a continuous linear
form on each
  $\Omega CS_{\rm cpt}^a(U)$ and is closed  on non integer order symbols
  valued forms:
$$\left( \alpha=d\, \beta  \in  \Omega CS^{\notin\Z}_{\rm cpt}(U)\right) \Rightarrow
\cutoffint_{T^*U}
  \alpha=0.$$
\end{lem}
{\bf Proof:} We prove the first statement. The continuity follows from the continuity  of the cut-off
regularised integral on $CS_{\rm c.c.}^a( \R^n)$ for any $a\in \C$. Similarly,  Stokes' property  follows from Stokes'
property of the ordinary   integral on $C_{\rm cpt}^{<-n}(U)$ combined with the fact
that the  cut-off regularised integral  $\cutoffint_{\R^n}$ vanishes on derivatives
$\partial_{\xi_j}$ of non integer order symbols as a result of
Proposition \ref{prop:Stokesint}. \endsquare\\ \\

Using a partition of unity, one can patch up the cut-off regularised integral
of symbols with compact support to a canonical trace on non integer order
classical pseudo-differential operators \cite{KV}.
\begin{defn}  The canonical trace  is  defined on $\Cl^{\notin \Z}(M, E)$ by
$${\rm TR}(A):=\frac{1}{(2\pi)^n}\, \int_Mdx \cutoffint_{T_x^*M}{\rm tr}_x\left(\sigma(A)(x, \xi)\right)\,
d\,\xi=\frac{1}{\sqrt{2\pi}^n}\,  \int_{M}{\rm TR}_x(A)\,d\, x  $$
where
${\rm TR}_x(A):=\frac{1}{\sqrt{2\pi}^n} \cutoffint_{T_x^*M }{\rm tr}_x\left(\sigma(A)(x, \xi)\right)\, d\xi$
 is the canonical trace density at point $x$.
\end{defn}
The canonical trace is tracial on non integer order operators as a
consequence of  Stokes' property for cut-off regularised integrals
on non integer order symbols.
\begin{prop} Let $A\in \Cl(M),\quad B\in \Cl(M,E)$ be two classical  operators with non integer order such that their bracket
  $[A, B]$ also has non integer order. Then
$${\rm TR}\left([A,B]\right)=0.$$
\end{prop}
{\bf Proof:}
The product of   $A$ and $ B$ on $M$
reads$$
A\, B=\sum_{\vert \alpha\vert \leq N}\frac{(-i)^{\vert \alpha\vert}}{\alpha!} {\rm Op}(\partial_\xi^\alpha \sigma(A)\,
\partial_x^\alpha \sigma(B))+ R_N(A\, B)$$
  for any integer $N$  and with   $R_N(A\, B)$ of order $<a+b-N$ where $a, b$
  are the orders of $A, B$ respectively.  Hence
$$
[A, B]=\sum_{\vert \alpha\vert \leq N}\frac{(-i)^{\vert \alpha\vert}}{\alpha!} {\rm Op}\left(\partial_\xi^\alpha \sigma(A)\,
\partial_x^\alpha \sigma(B)- \partial_\xi^\alpha \sigma(B)\,
\partial_x^\alpha \sigma(A)\right)+ R_N([A, B])$$
with similar notations.\\
When the bracket $[A, B]$ has non integer order, we can apply the canonical
trace on either side and write:
\begin{eqnarray*}
{\rm TR}\left([A, B]\right)&=&\sum_{\vert \alpha\vert \leq N}\frac{(-i)^{\vert
    \alpha\vert}}{\alpha!} \int_M d\,x\,
\cutoffint_{\R^n}{\rm tr}_x\left(\partial_\xi^\alpha
\sigma(A)\,
\partial_x^\alpha \sigma(B)- \partial_\xi^\alpha \sigma(B)\,
\partial_x^\alpha \sigma(A)\right)\, d\,\xi+ {\rm tr}\left(R_N([A,
B])\right)\\
&=& \sum_{\vert \alpha\vert \leq N}\frac{(-i)^{\vert
    \alpha\vert}}{\alpha!} \int_M d\,x\,
\cutoffint_{\R^n}{\rm tr}_x\left(\partial_\xi^\alpha
\sigma(A)\,
\partial_x^\alpha \sigma(B)-
\partial_x^\alpha \sigma(A)\, \partial_\xi^\alpha \sigma(B)\right)\, d\,\xi+ {\rm tr}\left(R_N([A,
B])\right)\\
&=& {\rm tr}\left(R_N([A,
B])\right).
\end{eqnarray*}
In the last identity,  we used Stokes' property for cut-off regularised integrals on non
integer order symbols (see Lemma \ref{lem:cutoffintStokesU}) to implement
repeated integration by parts in order to show
that   the integral
term on the r.h.s. vanishes using the fact that the ordinary trace on matrices
is cyclic. \\
Thus we have
$${\rm TR}\left([A,B]\right)={\rm tr}\left(R_N([A,
B])\right)$$
with $R_N([A,
B])$ of order $< a+b-N$.\\
Since $N$ can be chosen arbitrarily large, we have ${\rm TR}\left([A,B]\right)={\rm tr}\left(R_\infty([A,
B])\right)$ for some smoothing operator $R_\infty([A,
B])$.\\ On the other hand, for any smoothing operators $S, T$  the operators
$[S, B]$ and $[A, T]$ are smoothing and a direct check using the kernel
representation of these operators shows that
  ${\rm TR}\left([S, B]\right)= {\rm tr}\left([S, B]\right)=0$ and similarly, ${\rm
    TR}\left([A, T]\right)={\rm tr}\left([A, T]\right)=0.$
It follows that
$
{\rm TR}\left([A+S,
B+T]\right)={\rm TR}\left([A,
B]\right)$ leading to   $${\rm tr}\left(R_\infty([A+S,
B+T])\right)
={\rm tr}\left(R_\infty([A,
B])\right)$$
for any smoothing operators $S, T$. But this  means that the bilinear form
$(A, B)\mapsto {\rm tr}\left(R_\infty([A,
B])\right)$ is purely symbolic, namely that it  depends only on a finite
number of homogeneous components of the symbols of $A$ and $B$, which by its
very construction is clearly not the case unless it vanishes. This proves that
${\rm tr}\left(R_\infty([A,
B])\right)=0$ and hence that ${\rm TR}\left([A,
B]\right)=0$.\endsquare

\subsection{Holomorphic families of classical pseudodifferential operators}
The notion of holomorphic family of classical pseudodifferential operators first introduced by Guillemin in \cite{G1} and extensively
used by Kontsevich and Vishik in \cite{KV} 
generalises the notion of complex power $A^z$ of an elliptic operator developped by Seeley \cite{Se}, the
derivatives of which lead to logarithms.  
\begin{defn}
 Let $\Omega$ be a domain of $\C$ and $U$ an open subset of $\R^n$. A family $(\sigma(z))_{z\in \Omega}\subset CS(U)$ is holomorphic when \\
(i) the order $\alpha(z)$ of $\sigma(z)$ is holomorphic  on $\Omega$.\\
(ii) For $(x, \xi)\in U\times\R^n$, the function $z\to\sigma(z)(x,
\xi)$ is holomorphic on $\Omega$ and
 $\forall k \geq 0, \partial_z^k\sigma(z)\in S^{\alpha(z)+\e}(U)$ for all $\e>0$.\\
(iii) For any integer $j\geq 0,$ the (positively) homogeneous component
$\sigma_{\alpha(z)-j}(z)(x, \xi)$ of degree $\alpha(z)-j$ of the symbol  is holomorphic on $\Omega.$
\end{defn}
The derivative  of  a holomorphic family $\sigma(z)$ of classical
  symbols yields a holomorphic family of symbols, the  asymptotic expansions of
  which a priori involve a logarithmic term.
\begin{lem} The derivative of  a holomorphic family $\sigma(z)$ of classical
  symbols of order $\alpha(z)$ defines a
holomorphic family of symbols $\sigma^\prime(z)$  of order $\alpha(z)$ with
asymptotic expansion:
\begin{equation}\label{eq:sigmaprime}\sigma^\prime(z)(x, \xi)\sim\sum_{j=0}^\infty  \chi(\xi)\left( \log \vert
  \xi\vert\, \sigma_{\alpha(z)-j, 1}^\prime(z)(x, \xi)
+ \sigma^\prime_{\alpha(z)-j, 0}(z)(x, \xi)\right)\quad\forall (x, \xi)\in T^*U
\end{equation} for some smooth cut-off function $\chi$ around the origin which is
identically equal to $1$ outside the open unit ball and positively
homogeneous symbols
\begin{equation}\label{eq:sigmaprimej}\sigma_{\alpha(z)-j, 0}^\prime(z)(x, \xi)=\vert \xi\vert^{\alpha(z)-j}\,  \partial_z
\left(\sigma_{\alpha(z)-j}(z)(x,\frac{\xi}{\vert
  \xi \vert})\right), \quad 
\sigma^\prime_{{\alpha(z)-j}, 1}(z)=\alpha'(z)\, 
\sigma_{\alpha(z)-j}(z)(x,\xi)
\end{equation}
of degree $\alpha(z)-j$. 
\end{lem}{\bf Proof:} We write 
$$\sigma(z)(x, \xi)\sim\sum_{j=0}^\infty  \chi(\xi) \, \sigma_{\alpha(z)-j}(z)(x,
\xi).$$ Using  the positive homogeneity of the components
$\sigma_{\alpha(z)-j}$ we have:
\begin{eqnarray*}
&{}&\partial_z \left(\sigma_{\alpha(z)-j}(z)(x, \xi)\right)\\
&=& \partial_z \left(\vert \xi
\vert^{\alpha(z)-j}\sigma_{\alpha(z)-j}(z)(x,\frac{\xi}{\vert \xi \vert}) \right)\\
&=& \left(\alpha'(z)\vert \xi \vert^{\alpha(z)-j}\sigma_{\alpha(z)-j}(z)(x,\frac{\xi}{\vert \xi
\vert}) \right)\log \vert \xi \vert
 +\vert \xi
\vert^{\alpha(z)-j} \partial_z
\left(\sigma_{\alpha(z)-j}(z)(x,\frac{\xi}{\vert
  \xi \vert})\right)\\
&=& \left(\alpha'(z)\sigma_{\alpha(z)-j}(z)(x,\xi) \right)\log \vert \xi \vert
 +\vert \xi
\vert^{\alpha(z)-j} \partial_z
\left(\sigma_{\alpha(z)-j}(z)(x,\frac{\xi}{\vert
  \xi \vert})\right)
\end{eqnarray*}
which shows that $\partial_z \left(\sigma_{\alpha(z)-j}(z)(x, \xi)\right)$ has
order $\alpha(z)-j$. 
Thus 
$$\partial_z\left( \sigma_N(z)(x, \xi)\right)= \sigma^\prime(z) (x, \xi)- \sum_{j<N} \chi(\xi)\,
\partial_z\left(\sigma_{\alpha(z)-j}(z)(x,  \xi)\right)$$
lies in ${\cal S}^{\alpha(z)-N+\e}(U)$ for any $\e>0$ so that
$\sigma^\prime(z)$ is a symbol of order $\alpha(z)$ with asymptotic expansion:
\begin{equation}\sigma^\prime(z)(x, \xi)\sim
  \sum_{j=0}^\infty \chi(\xi)\, \sigma_{\alpha(z)-j}^\prime(z) \quad\forall (x, \xi)\in
  T^*U\end{equation}
where
$$\sigma^\prime_{\alpha(z)-j}(z)(x, \xi):=\log \vert \xi\vert
 \sigma^\prime_{\alpha(z)-j,1}(z)(x, \xi)+  \sigma^\prime_{\alpha(z)-j,0}(z)(x, \xi)$$
for some positively homogeneous symbols $$\sigma^\prime_{\alpha(z)-j,0}(z)(x,
\xi):=\vert \xi\vert^{\alpha(z)-j}\,  \partial_z
\left(\sigma_{\alpha(z)-j}(z)(x,\frac{\xi}{\vert
  \xi \vert})\right)$$
and $$\sigma^\prime_{\alpha(z)-j,1}(z)(x, \xi):= \alpha^\prime(z)\,\sigma_{\alpha(z)-j}(z)(x, \xi) $$ of degree $\alpha(z)-j$. \\
 On the other hand, differentiating the asymptotic expansion $\sigma(z)(x, \xi)\sim \sum_{j=0}
\chi(\xi)\, \sigma_{\alpha(z)-j} (z)(x, \xi) $  w.r.
to $z$ yields
 $$
\sigma^\prime(z)(x, \xi)\sim \sum_{j=0} \chi(\xi)\,
\partial_z\left(\sigma_{\alpha(z)-j} (z)(x, \xi)\right). $$
Hence, 
$$\partial_z\left(\sigma_{\alpha(z)-j}(z) (x,
  \xi)\right)=\sigma^\prime_{\alpha(z)-j}(z)(x, \xi)=\vert \xi\vert^{\alpha(z)-j}\,  \partial_z
\left(\sigma_{\alpha(z)-j}(z)(x,\frac{\xi}{\vert
  \xi \vert})\right)+   \alpha^\prime(z)\,\sigma_{\alpha(z)-j}(x, \xi)\, \log
\vert \xi\vert $$ as announced.
\endsquare\\ \\
Correspondingly we recall the notion of holomorphic classical
pseudodifferential operators.
\begin{defn} A family $(A(z))_{z \in \Omega}\in \Cl(M, E)$ is holomorphic if
  in any local trivialisation   we can write $A(z)$ in the form
 $A(z)=Op(\sigma(A(z)))+R(z)$, for some holomorphic family of symbols
 $\left(\sigma(A(z))\right)_{z\in \Omega}$ and some holomorphic family $(R(z))_{z \in
   \Omega}$ of smoothing operators i.e. given by a holomorphic family of
 smooth Schwartz kernels.
\end{defn}
It follows from (\ref{eq:sigmaprime})  and (\ref{eq:sigmaprimej})  that
\begin{eqnarray}\label{eq:Aprimej}&{}&\partial_z\left(\sigma(A(z))_{\alpha(z)-j}\right)(x,
  \xi)=\sigma_{\alpha(z)-j}(A^\prime(z))(x, \xi)\nonumber
  \\
&=&\alpha^\prime(z)
  \sigma_{\alpha(z)-j}(A(z))(x, \xi)\, \log\vert \xi\vert+ \vert \xi\vert^{\alpha(z)-j}\,  \partial_z
\left(\sigma_{\alpha(z)-j}(A(z))(x,\frac{\xi}{\vert
  \xi \vert})\right)(x, \xi).
\end{eqnarray}
We call admissible with spectral cut $\theta$ an operator
$A\in\Cl(M, E)$ with  leading symbol $\sigma_L(A)$ that has no
eigenvalue on the ray $L_\theta=\{re^{i\theta}, r\geq 0\}$ in which
case it is elliptic, and such that the spectrum of $A$ does not meet
the open ray $\{re^{i\theta}, r> 0\}$. In that case,
 following Seeley \cite{Se}, one can define the complex power $A_\theta^z$
 which yields a holomorphic family $z\mapsto A_\theta^z$  in $\Cl(M,E)$.
\begin{ex} Given an admissible operator $A\in \Cl^a(M, E)$ with spectral cut
  $\theta$,  the operator
$A(z)=A_\theta^{z}$ is a holomorphic family and we have
 $A^\prime(0)=(\partial_z A_\theta^{z})_{z=0}=\log_{\theta}A $. Furthermore,
 it follows from  
(\ref{eq:Aprimej}) that: 
$$\sigma_{-j}(\log_\theta A)(x, \xi)= \partial_z\left(
  \sigma_{az-j}(A_\theta^z)\right)_{\vert_{z=0}}(x, \xi) = a \, \delta_{j, 0}\,
\log\vert \xi\vert +
\sigma_{-j,0}(\log_\theta A) (x, \xi)$$
with $$\sigma_{-j, 0}(\log_\theta A)(x, \xi)= \vert \xi\vert^{-j}\,  \left(\partial_z
\left(\sigma_{az-j}(A_\theta^z)(x,\frac{\xi}{\vert
  \xi \vert})\right)(x, \xi)\right)_{\vert_{z=0}}.$$
\end{ex}
\subsection{Continuity of the canonical trace on non integer order classical
  pseudodifferential operators} It follows from the very definition of the canonical trace that it is
continuous  w.r. to the Fr\'echet
topology on $\Cl^a(M, E)$ for every $a\notin\Z$. In this paragraph we discuss its continuity on
(holomorphic) families of varying order.\\
The following proposition collects results from \cite{KV} and   \cite{PS}.
 \begin{prop}\label{prop:PS} Let $\sigma(z)\in CS(U)$ (resp. $A(z)\in \Cl(M, E)$)
   be a holomorphic family of order $\alpha(z)$ such that
   $\alpha^\prime(0)\neq 0$. The map 
$$z\mapsto \int_{T_x^*U} \sigma(z)(x, \xi)\, d\xi$$
(resp. $z\mapsto {\rm tr}(A(z))$)
is holomorphic on $\alpha^{-1}\left(]-\infty, -n[\right)$
and extends to a meromorphic map  $z\mapsto  \cutoffint_{T_x^*U}{\rm tr}_x\left(\sigma(z)(x, \xi)\right)\,
d\xi$ (resp. $z\mapsto {\rm TR}(A(z))$) on the complex plane with
simple poles and \cite{KV} \begin{equation}\label{eq:KVsymb}{\rm
       Res}_{z=0} \cutoffint_{Tx^*U}\sigma(z)(x, \xi)\, d\xi= -
  \frac{1}{\alpha^\prime(0)} {\rm res}_x(\sigma(0)(x, \xi)),
\end{equation} 
$($ resp. \begin{equation}\label{eq:KVop}{\rm Res}_{z=0} {\rm TR}(A(z))= -
  \frac{1}{\alpha^\prime(0)} {\rm res}(A(0)).\quad)
\end{equation}
 Furthermore, if $\alpha(z)$ is affine in $z$ \cite{PS} \begin{equation}\label{eq:PSsymb}{\rm
       fp}_{z=0} \cutoffint_{T_x^*U}\sigma(z)(x, \xi)\, d\xi= \cutoffint_{T_x^*U}{\rm tr}_x\left(\sigma(0)(x, \xi)\right)\, d\xi-
  \frac{1}{\alpha^\prime(0)} \int_{S_x^*U}\sigma^\prime(0)(x, \cdot)\, d\xi \quad\forall
x\in U
\end{equation}
$($ resp.
\begin{equation}\label{eq:PSop}{\rm fp}_{z=0} {\rm TR}(A(z))= \int_M dx\left({\rm  TR}_x(A(0))-
  \frac{1}{\alpha^\prime(0)} {\rm res}_x(A^\prime(0))\right).\quad)
\end{equation}
\end{prop}
\begin{cor} The canonical trace on non integer order operators is continuous
  along  holomorphic families  with
  affine order. In other words, for any holomorphic family $A(z)\in \Cl(M, E)$
  with  affine order
  $\alpha(z)$ such that $A(0)\in \Cl^{\notin \Z}(M, E)$
$$\lim_{z\to 0} {\rm TR}(A(z))= {\rm TR}(A(0)).$$
\end{cor}

{\bf Proof:} We can assume that the order $\alpha(z)$ of $A(z)$ satisfies
 $\alpha^\prime(0)\neq 0$ for otherwise the order is constant in which case we
 already know that the canonical trace is continuous at $0$.\\
 If $\alpha^\prime(0)\neq 0$, the map  $z\mapsto {\rm TR}(A(z))$ is holomorphic at $z=0$  since
 by equation (\ref{eq:KVop})$${\rm Res}_{z=0}{\rm TR}(A(z))= \frac{1}{\alpha^\prime(0)} {\rm
  res}(A(0))$$  which vanishes as a result of the non integrality of the order
of  $A(0)$. Since the  derivative
$A^\prime(0)$ at $z=0$ has same   order $\alpha(0)$ as $A(0)$
which is non integer by assumption, $A^\prime(0)$ also  has vanishing residue
density so that by equation  (\ref{eq:PSop}) we have:
$${\rm fp}_{z=0} {\rm TR}(A(z))= \frac{1}{\sqrt{2\pi}^n}\int_M  {\rm TR}_x(A(0))\,dx= {\rm TR}(A(0)).$$
\endsquare
\subsection{Odd- (resp. even-) class operators embedded in holomorphic families}
 For any integer $a$, the condition
  $ \sigma_{a-j}(x, -\xi)= (-1)^{a-j}\sigma_{a-j}(x,
  \xi)\quad \forall j\in \N\cup \{0\}$ which characterises a classical symbol
  of order $a$ that lies  in the
  odd-class, extends to
  log-polyhomogeneous symbols of logarithmic type $1$: $$\sigma(x, \xi)\sim \sum_{j=0}\chi(\xi)\,
  \sigma_{a-j}(x, \xi); \quad \sigma_{a-j}= \sigma_{a-j}^0+\sigma_{a-j}^1\, \log \vert
  \xi\vert$$ with $\sigma_{a-j}^i(x, \cdot), i=0, 1$ positively homogeneous of degree $a-j$. One   requires that
  $$\sigma_{a-j}(x, -\xi)=(-1)^{a-j} \sigma_{a-j}(x, \xi)\quad \forall j\in \N\cup
  \{0\}\quad \forall (x, \xi)\in T^*U$$ or equivalently that
  $$\sigma^i_{a-j}(x, -\xi)= (-1)^{a-j}\sigma^i_{a-j}(x,
  \xi) \quad \forall j\in \N\cup
  \{0\}\quad \forall (x, \xi)\in T^*U$$ for both
  $i=0$ and $i= 1$.
\begin{prop}\label{prop:Braverman}\cite{B} Given an  admissible operator  $A\in
  \Cl^{\rm odd}(M,E)$  with positive order $a>0$ and spectral cuts
  $\theta$ and $\theta-a\pi$, the symmetrised logarithm
  $$A^\prime_\theta(0)= \frac{\log_{\theta}A+\log_{\theta-a\pi}A}{2}$$  where  we
have set $A_\theta(z) :=\frac{ A_\theta^z+
    A_{\theta-a\pi}^z}{2}$,
 lies in the odd-class.

\end{prop}
\begin{rk} When the order $a$ of $A$ is even, then $A_\theta= A_{\theta-a\pi}$
  so that  $A_\theta(z)=A_\theta^z$ and $A_\theta^\prime(0)=\log_\theta
  A$. This yields back the known fact \cite{KV} that the logarithm of an
  odd-class admissible \pdo with even order lies in the odd-class. 
\end{rk}
{\bf Proof:} Recall that
 the homogeneous components of the symbol of
$A_\theta^{z}$ are 
\begin{equation}\label{eq:homComponentComplexP}
\sigma_{az-j}(A_\theta^{z})(x,\xi)=\frac{i}{2\pi}
\int_{\Gamma_\theta} \lambda_\theta^{z}\, q_{-a-j}(x,\xi,\lambda)\,
d\lambda.
\end{equation}
with\cite{S}
$$q_{-a}=(\sigma_a(A)-\lambda)^{-1}$$
$$q_{-a-j}=-q_{-a}\left(\sum_{k+l+\vert
\alpha \vert =j,l<j}
\frac{1}{\alpha!}\partial_{\xi}^{\alpha}\sigma_{a-k}(A)
D_x^{\alpha}q_{-a-l}\right)$$
 the positively homogeneous components of the resolvent $(A-\lambda I)^{-1}$.
In other words, these components $q_{-a-j}$ are  positively
homogeneous   in $(\xi,\lambda^{\frac{1}{a}})$  i.e. for $t>0,$ for
$(x, \xi)\in T^\star M,$ 
\begin{equation}\label{eq:homComponentResolvent}
q_k(x,t\xi,t^{\frac{1}{a}}\lambda)=t^kq_k(x,\xi,\lambda)\quad\forall t>0.
\end{equation}
 If $A\in\Cl^a(M,E)$ lies in the odd-class,  this extends to  any real
 number $t$ since  we have \cite{KV} par. 2
  \begin{equation}\label{eq:oddClassHomResolvent}
q_k(x,-\xi,(-1)^{a}\lambda)=(-1)^k\, q_k(x,\xi,\lambda).
\end{equation}
Now, assume that ${\rm Re}\,z<0$. A Cauchy integral gives
\begin{eqnarray*}
\sigma_{za-j}(A_\theta^z)(x, -\xi) &=&\frac{i}{2\pi}
\int_{\Gamma_\theta} \lambda_\theta^{z}\,q_{-a-j}(x,-\xi,\lambda)\,
d\lambda\\
&=&(-1)^{a+j}\frac{i}{2\pi} \int_{\Gamma_\theta}
\lambda_\theta^{z}\,q_{-a-j}(x,\xi,(-1)^a\lambda)\,
d\lambda\\
\end{eqnarray*}
where $\Gamma_\theta$ is an appropriate contour around the spectrum
of $A$ along
the ray $L_\theta$.\\
By a change of variable, we obtain
\begin{eqnarray*}
\sigma_{az-j}(A_\theta^z)(x, -\xi) &=&(-1)^{a+j}\frac{i}{2\pi}
\int_{\Gamma_{\theta-a\pi}}
(e^{ia\pi}\mu)_{\theta}^{z}\,q_{-a-j}(x,\xi,\mu)\,
d(e^{ia\pi}\mu)\\
&=&(-1)^{a+j}e^{ia\pi}\frac{i}{2\pi} \int_{\Gamma_{\theta-a\pi}}
e^{iza\pi}\mu_{\theta-a\pi}^{z}\,q_{-a-j}(x,\xi,\mu)\,
d\mu\\
&=&(-1)^je^{iaz\pi}\sigma_{az-j}(A_{\theta-a\pi}^z)(x, \xi).
\end{eqnarray*}
Thus
\begin{equation}\label{eq:B1}\sigma_{az-j}(A_\theta^z)(x, -\xi)
=e^{i\pi(az-j)}\sigma_{az-j}(A_{\theta-a\pi}^z)(x, \xi).
\end{equation}
Since
both the left and the right hand side of this equality are analytic
in $z$, we conclude that the equality holds for all $z\in \C.$
Similarly, \begin{eqnarray*}
\sigma_{az-j}(A_{\theta-a\pi}^z)(x, -\xi) &=&(-1)^{a+j}\frac{i}{2\pi}
\int_{\Gamma_{\theta-a\pi}}
(e^{-ia\pi}\mu)_{\theta-a\pi}^{z}\,q_{-a-j}(x,\xi,\mu)\,
d(e^{-ia\pi}\mu)\\
&=&(-1)^{a+j}e^{-ia\pi}\frac{i}{2\pi} \int_{\Gamma_{\theta}}
e^{-iza\pi}\mu_\theta^{z}\,q_{-a-j}(x,\xi,\mu)\,
d\mu\\
&=&(-1)^je^{-iaz\pi}\sigma_{az-j}(A_{\theta}^z)(x, \xi)
\end{eqnarray*}
so that
\begin{equation}\label{eq:B2}\sigma_{az-j}(A_{\theta-a\pi}^z)(x, -\xi)
=e^{i\pi(-az-j)}\sigma_{az-j}(A_{\theta}^z)(x, \xi).
\end{equation}
\\
Differentiating (\ref{eq:B1})  w.r. to $z$ on either side yields:
\begin{eqnarray}\label{eq:B3}&{}&\sigma_{az-j}(\partial_z A_\theta^z)(x, -\xi)\\
&=&\partial_z \left(\sigma_{az-j}(A_\theta^z)(x, -\xi)\right)\nonumber\\
&=& \partial_z \left(e^{i\pi(az-j)}\sigma_{az-j}(A_{\theta-a\pi}^z)(x,
  \xi)\right)\nonumber\\
&=& i\pi \, a  \,e^{i\pi(az-j)}\,\sigma_{az-j}(A_{\theta-a\pi}^z)(x,
  \xi)+e^{i\pi(az-j)}\sigma_{az-j}(\partial_z A_{\theta-a\pi}^z)(x,
  \xi).
\end{eqnarray}
Similarly, differentiating (\ref{eq:B2}) on either side yields:
\begin{eqnarray}\label{eq:B4}&{}&\sigma_{az-j}(\partial_z A_{\theta-a\pi}^z)(x, -\xi)\\
&=&\partial_z \left(\sigma_{az-j}(A_{\theta-a\pi}^z)(x, -\xi)\right)\nonumber\\
&=& \partial_z \left(e^{i\pi(-az-j)}\sigma_{az-j}(A_{\theta}^z)(x,
  \xi)\right)\nonumber\\
&=& -i\pi \, a  \,e^{i\pi(az-j)}\, \sigma_{-az-j}(A_{\theta}^z)(x,
  \xi)+e^{i\pi(-az-j)}\sigma_{-az-j}(\partial_zA_{\theta}^z)(x,
  \xi).
\end{eqnarray}
Combining equations (\ref{eq:B3}) and (\ref{eq:B4}) yields at $z=0$:
\begin{eqnarray*}
\sigma_{-j}(A^\prime(0))(x, -\xi)&=& \frac{
\sigma_{-j}\left(\left( \partial_z A_\theta^z\right)_{\vert_{z=0}}\right)(x,-
\xi)+
\sigma_{-j}\left(\left( \partial_z
    A_{\theta-a\pi}^z\right)_{\vert_{z=0}}\right)(x, -\xi)}{2}\\
&=&\frac{  i\pi \, a  \,\delta_{j,0}+(-1)^j\sigma_{-j}\left(\left(\partial_z A_{\theta-a\pi}^z\right)_{\vert_{z=0}}\right)(x,
  \xi)- i\pi \, a  \,\delta_{j,0}+(-1)^j\sigma_{-j}\left(\left(\partial_z
      A_{\theta}^z\right)_{\vert_{z=0}}\right)(x,
  \xi)}{2}\\
&=& (-1)^j\sigma_{-j}(A^\prime(0))(x, \xi)
\end{eqnarray*}
so that $A^\prime(0)$ lies in the odd-class. 
\endsquare
\begin{ex} Take $M$ a Riemannian manifold and $A= \Delta_g$ the Laplace
  Beltrami operator with $\theta=\frac{\pi}{2}$. It has order $2$ and lies in the odd-class
  since it is a differential operator.  Then $A_\theta^\prime(0)= 
   \log_{\frac{\pi}{2}}  \Delta_g$ lies in the odd-class.
\end{ex}
\begin{ex}  Let $M$ be a spin manifold and  $E=S\otimes W$ a twisted
bundle with $S$ the spinor bundle and $W$ an exterior vector bundle
over  $M$. From a twisted connection $\nabla^E=
\nabla^S\otimes 1+ 1\otimes \nabla^W$ on $E$, with $\nabla^S$  the
connection  on $S$ induced by the Levi-Civita connection on $M$,
$\nabla^W$ a connection on $W$, one can build the corresponding
twisted Dirac operator $D= c\cdot \nabla^E$. Here $c$ stands for the
Clifford multiplication on the Clifford bundle $E$. $D$  is an
admissible operator of order $1$ which lies in the odd-class since
it is a differential operator. Moreover, since it is self-adjoint,
its spectrum is real so that it has spectral cut
$\theta:=\frac{\pi}{2}$ and $\theta-\pi= -\frac{\pi}{2}$.\\
Take the order one differential operator $A= D$. It lies in the odd-class and
so does its symmetrised logarithm $\frac{\log_{\frac{\pi}{2}}D+ \log_{-\frac{\pi}{2}}D}{2}$.
\end{ex}
\begin{cor}\label{cor:Braverman} Provided there is
    an admissible operator $Q\in \Cl^{\rm odd}(M, E)$ with positive
    order $q$ and spectral cuts $\theta$ and $\theta-q\pi$,then
    any operator $A\in \Cl^{\rm odd} (M, E)$ (resp. $A\in \Cl^{\rm even} (M, E)$) can be embedded in a
    holomorphic family $$A_\theta^Q(z):= A\, \frac{Q_\theta^z+
    Q_{\theta-q\pi}^z}{2}$$ such that $\left(\partial_z
    A_\theta^Q(z)\right)_{\vert_{z=0}}$ lies in the odd-class (resp. even-class).
\end{cor}
{\bf Proof:} This follows from  Proposition \ref{prop:Braverman} applied to
$Q$ combined with the stability of the odd-class under
products (resp. the fact that the product of an  even and odd-class operator
is even).\\ Let us focus on the odd-class case, since the proof in the
even-class  case goes in
a similar manner.  
\\ By Proposition \ref{prop:Braverman}, $Q_\theta^\prime(0)$ lies in
the odd-class. Since $$\left(\partial_z
    A_\theta^Q(z)\right)_{\vert_{z=0}}= A\,       Q_\theta^\prime(0), $$
 applying 
$$\sigma(AB)=\sum_{\alpha }\frac{(-i)^\alpha}{\alpha!} \partial^\alpha_\xi\sigma(A) \partial_x^\alpha
\sigma(B)$$ to $A$ and $B= Q_\theta^\prime(0)$ yields that $\left(\partial_z
    A_\theta^Q(z)\right)_{\vert_{z=0}}$
lies in the odd-class since $A\in \Cl^{\rm odd}(M, E)$. 
\endsquare
\subsection{The canonical trace  on  odd (resp. even)-class operators
  in odd (resp. even) dimensions}
{\it In the sequel, $M$ is an
odd-(resp. even-) dimensional  manifold. Let $\pi: E\to M$ be a vector bundle over $M$ such
that there is
    an admissible operator $Q\in \Cl^{\rm odd}(M, E)$ with positive
    order $q$ and spectral cuts $\theta$ and $\theta-q\pi$. }  
\begin{rk}  In view of the above examples, these assumptions are fulfilled in very natural geometric setups. 
\end{rk}
\begin{thm}\label{thm:context} The canonical trace ${\rm  TR}$  extends continuously to
  $\Cl^{\rm odd}(M, E)$ in odd dimensions (resp.  $\Cl^{\rm even}(M, E)$ in
  even dimensions)  in the following manner. \\ Let $M$ be odd (resp. even)
  dimensional. For any  holomorphic family
  $A(z)\in \Cl(M, E)$   with non
  constant affine order such that both   $A(0)$ and $ A^\prime(0)$ lie in
  $\Cl^{\rm odd}(M,E)$ (resp.  $\Cl^{\rm even}(M,E)$),
\begin{enumerate}
\item  the map $z\mapsto {\rm
  TR}(A(z))$ is holomorphic at $z=0$,
\item $\int_{M}\left(\int_{T_x^*M}{\rm
  tr}_x\sigma(A(0))(x, \xi)\, d\xi\right)\, dx$ defines a global density on $M$
so that  $${\rm TR}(A(0))=\frac{1}{(2\pi)^n}\, \int_{M}\left(\int_{T_x^*M}{\rm
  tr}_x\sigma(A(0))(x, \xi)\, d\xi\right)\, dx
$$ is well-defined, \item 
$\lim_{z=0} {\rm TR}(A(z))= {\rm TR}(A(0)).$
\end{enumerate} 
\end{thm}
{\bf Proof: } We carry out the proof in odd dimensions for odd-class
operators. The proof goes similarly in the even dimensional case for
even-class operators.
\\
Since the noncommutative residue vanishes on  $\Cl^{\rm
  odd}(M,E)$ and $A(0)\in \Cl^{\rm
  odd}(M,E)$, we have ${\rm res}(A(0))=0$.  It follows from (\ref{eq:KVsymb}) that the complex residue ${\rm Res}_{z=0}{\rm
  TR}(A(z)) $ which is proportional to the noncommutative reside res$(A(0))$
vanishes so that the map $z\mapsto {\rm
  TR}(A(z)) $ is holomorphic at $z=0$. We now apply (\ref{eq:PSsymb})  to
$\sigma(z):=\sigma(A(z))$; since $A^\prime(0)$ lies in the odd-class, it has
vanishing residue density ${\rm res}_x(A^\prime(0))$. Consequently,  $$\lim_{z=0} \cutoffint_{T_x^*M}{\rm tr}_x(\sigma(A(z))(x, \xi))\, d\xi=
     \cutoffint_{T_x^*M}{\rm tr}_x(\sigma(A)(x, \xi))\, d\xi= \sqrt{2\pi}^n\,{\rm TR}_x(A(0))\quad\forall x\in
     M.$$ Since the l.h.s gives rise
     to  a
     globally defined  density $$\left({\rm
       fp}_{z=0} \cutoffint_{T_x^*M}\sigma(A(z))(x, \xi)\, d\xi\right)\,
   dx=\left(\lim_{z\to 0} \cutoffint_{T_x^*M}\sigma(A(z))(x, \xi)\, d\xi\right)\, dx$$
   so does the right hand side give rise to a globally defined density  ${\rm
     TR}_x(A(0))\, dx$. Integrating over $M$ yields the existence of ${\rm
     TR}(A(0))$ and:
$$\lim_{z\to 0} {\rm TR}(A(z))=\frac{1}{\sqrt{2\pi}^n} \int_M
{\rm  TR}_x(A(0))\, dx= {\rm TR}(A(0)).$$\endsquare
\\ \\
\begin{cor}Any operator  $A\in \Cl^{\rm odd}(M, E)$ in odd dimensions (resp. $A\in
  \Cl^{\rm even}(M, E)$ in even dimensions)  has well-defined canonical trace
  $${\rm TR}(A)=\int_M dx\left( \cutoffint_{T_x^*M } {\rm tr}_x\left(
      \sigma_A(x, \xi)\right)\right)\, d\xi$$
and 
${\rm TR}(A)= \lim_{z\to 0} {\rm TR}(A(z))$ for any holomorphic family $A(z)$
with non constant affine order  such that $A(0)=A$ and $ A^\prime(0)$ lie in $\Cl^{\rm
  odd}(M, E)$ (resp. $\Cl^{\rm
  even}(M, E)$). In particular, 
\begin{itemize}
\item Kontsevich and Vishik's (resp. Grubb's) extended canonical trace
  \cite{KV}  (resp. \cite{Gr}) on odd-class (resp. even-class)
  operators in odd (resp. even) dimensions,
  $$A\mapsto {\rm Tr}_{(-1)}(A):=
  \lim_{z\to 0} {\rm TR}(A\, Q_{\frac{\pi}{2}}^z)$$
with  $Q\in \Cl^{\rm odd}(M, E)$ an admissible operator  of even positive order close enough
to a positive self-adjoint operator,
\item the symmetrised
    trace introduced by Braverman \cite{B} on odd-class operators in odd dimensions
$$A\mapsto {\rm Tr}^{\rm sym}(A):= \lim_{z\to 0}{\rm TR}(A_\theta^Q(z))$$
with  $Q\in \Cl^{\rm odd}(M, E)$ an admissible operator  of any positive order
$q$ 
 and with spectral cuts $\theta$ and $\theta-q\pi$,

\end{itemize}
 coincide with the
    canonical trace.
\end{cor}
We also recover  as a side result the  fact \cite{KV}
 that the extended canonical trace  vanishes on
brackets of odd-class operators.
\begin{cor} In odd dimensions, for any operators  $A\in \Cl^{\rm odd}(M, E)$,  $B\in\Cl^{\rm odd}(M, E)$  we have
$${\rm TR}\left([A, B]\right)=0.$$
\end{cor}
{\bf Proof:} With the notations of Proposition \ref{prop:Braverman},
the holomorphic family $C(z)=\left[A_\theta^Q(z)  ,
B_\theta^Q(z)\right]$ has derivative $C^\prime(0)=\left[A\, Q_\theta^\prime(0)  ,
B\right]+\left[A,B \,  Q_\theta^\prime(0)\right] $ at $z=0$. It  lies in the
odd class as a result of the stability of the odd-class under products. The result then follows
from applying Theorem \ref{thm:context}  to the holomorphic family
$C(z)$. Since ${\rm TR}\left(C(z)\right)$ vanishes as a meromorphic
map as a consequence of the vanishing of the canonical trace on non
integer order brackets, taking finite parts as $z\to 0$ we find that
$$0=  {\rm TR}([A, B]).$$
\endsquare
\vfill \eject\noindent
\section{Uniqueness : Characterisation  of linear forms that vanish on 
operator brackets}
We prove that the noncommutative residue on the algebra of 
classical pseudodifferential operators, the canonical trace on the set of  non integer
order ones, or in odd (resp. even) dimensions on the classes of
odd-(resp. even-) class
operators, are the unique (possibly continuous) linear forms that vanish on
brackets. The classes to which the canonical trace naturally extends have in
common that their operators have symbols with vanishing residue density. 
\subsection{Uniqueness of  the noncommutative residue }
We use the notations of section 2; in particular $U$ is an open connected
subset of $\R^n$.
\begin{prop}\label{prop:uniquenessresU} Any singular linear form $\rho:
  CS_{\rm cpt}(U)\to \C$ which restricts to
a continuous map on  $CS_{\rm cpt}^a(U)$  for any a $\in \C$ and which
fulfills Stokes' property  is proportional to the noncommutative residue.\\
Equivalently, any closed singular linear form $\tilde \rho:
 \Omega CS_{\rm cpt}(U)\to \C$ which restricts to
a continuous linear form on  $\Omega CS_{\rm cpt}^a(U)$ for any a $\in \C$  is
proportional to the  noncommutative residue $\widetilde {\rm res}$ extended to forms.
\end{prop}
{\bf Proof:} By similar arguments as in the case of symbols with
constant coefficients we can check that the  two statements are
equivalent. Let us  prove the first one. 
For any fixed $f\in C^\infty_{\rm cpt}(U)$ the map
$\rho_f:\sigma\mapsto  \rho(f\, \sigma)$ defines a singular linear form on $CS_{\rm c.c}(\R^n)$
which vanishes on derivatives in $\xi$ since we have $\rho(f\, \partial_{\xi_j} \sigma)=  \rho(\partial_{\xi_j}(f\,  \sigma))=0$. By Theorem \ref{thm:uniquenessres}, it
follows that  there is a constant $c(f)$   such that $\rho(f\,
\sigma)=c(f)\, {\rm res}(\sigma)$ for any $\sigma\in CS_{\rm c.c}(\R^n)$. Since $f\mapsto
\rho(f\, \sigma)$ is continuous, $c:f\mapsto c(f)$ lies in
$\left(C^\infty_{\rm cpt}(U)\right)^\prime$.
A general symbol
 $\sigma\in  CS_{\rm cpt}(U)$ can be approximated by linear combinations of tensor products
 $f\otimes \sigma$ with $f\in C_{\rm cpt}^\infty(U)$, $\sigma \in CS_{\rm c.c}(\R^n)$. It follows  from the continuity of $\rho$ that there is a  distribution
 $F \in \left( C^\infty_{\rm cpt}(U)\right)^\prime$ 
 such that   $\rho( \sigma)=F( {\rm
res}(\sigma(x,\cdot))$ for any
$\sigma\in CS_{\rm cpt}(U)$. This distribution being continuous, it reads 
$F(f)= \int_U\psi(x) \, f(x)\, dx$ for some $\psi\in C^\infty(U)$ so that 
$$\rho(\sigma)= \int_U \psi(x)\, {\rm res}(\sigma(x, \cdot))\, dx.$$ But since
$\rho$ is closed by assumption, for any $\sigma=f\otimes \tau$ with $\tau  \in
CS_{\rm c.c}(\R^n)$ and 
$f\in \Ci(U)$
we have  $$0= \rho\left(\partial_{x_i} \left(f\otimes
    \tau\right)\right)=\rho\left(\partial_{x_i} f\otimes
  \tau\right).$$ Choosing $\tau$ with non vanishing residue and 
 integrating by parts implies that 
$$ \int_U\partial_{x_i} \left(\psi(x)\,f(x)\right)\, dx =0\quad \forall f\in C^\infty_{\rm
  cpt}(U)$$ so that $\psi$ is a constant $c$ and $\rho(\sigma)=c\, \int_U
{\rm res}(\sigma(x, \cdot))\, dx$ is proportional to the noncommutative residue.
\endsquare

We now derive from the characterisation of the residue on symbols in terms of
 Stokes' property, the   uniqueness (up to a  multiplicative constant) of the residue as
a trace on $\Cl(M)$ which restricts to continuous linear forms on each
$\Cl^a(M)$. It uses the following lemma.\\
\begin{lem}\label{lem:smoothing}(\cite{Po1} Lemma 3.20 and \cite{Po2} Lemma 4.4.) Any smoothing operator $A\in
  \Cl(M)$  can be written as a finite sum
of brackets $\sigma=\sum_{i=1}^n [x_i, B_i]$ with $B_i\in \Cl^{-n+1}(M)$.
\end{lem} {\bf Proof:}  We briefly sketch  the proof which we take from  \cite{Po1}  and \cite{Po2}. A smoothing operator $R$ has smooth kernel $k_R(x, y)$ so that $k_R(x,
y)-k_R(x, x)$ is smooth and vanishes on the diagonal. It follows that there
are smooth functions $k_{1}, \cdots, k_{n}$ such that
$k_R(x, y)= k_R(x, x)+ \sum_{j=1}^n (x_j-y_j) \, k_{j}(x, y).$  Let $Q$ be
the operator defined by the kernel $k_Q(x,y)= k_R(x, x)$  and let
$R_j,j=1,\cdots, n$ be the smoothing
operators defined by the kernels $k_j(x, y)$ then
$R= Q+\sum_{j=1}^n [x_j, R_j]$.\\  Set $H_j(x, y):=
y_j\vert y\vert^{-2} \, k_Q(x, x)$ and  let  $Q_j$ be the operator with kernel
$(x, y)\mapsto  H_j(x, x-y) $; by propsotion 2.7 in \cite{Po2}, it is a classical pseudodifferential operator of order $-n+1$.  Since
$$\sum_{j=1}^n (x_j-y_j) \, H_j(x,x- y)= \sum_{j=1}^n \frac{(x_j-y_j)^2}{\vert
  x-y\vert^2}\, k_R(x,x) = k_Q(x,y),$$ it follows that $Q=\sum_{j=1}^n [x_j,
Q_j]$. \\ Since $R_j$ are smoothing and
$Q_j$ are of order $-n+1$ the result of the lemma follows. \endsquare

\begin{thm}\label{thm:uniquenessresM} Any linear form $L:\Cl(M)\to \C$ which
  restricts to continuous linear forms on $\Cl^a(M)$ for any $a\in \C$ and
  which  vanishes on brackets
$$L\left([A, B]\right)=0\quad
  \forall A, B\in \Cl(M)$$ is proportional to the noncommutative residue.
\end{thm}{\bf Proof:}  By Lemma \ref{lem:smoothing}, such a linear form $L$
vanishes on smoothing operators. \\
 Given  a    local chart $(U,
\phi)$ on $M$,  the  map $$\rho_\phi:=L\circ\phi^*\circ {\rm
  Op}$$ then defines a singular linear form on   $CS_{\rm cpt}(\phi(U))$. \\
For any $\sigma \in CS_{\rm
  cpt}(\phi(U))$ and for any
  $x_j, j=1, \cdots, n$ corresponding to  the coordinates in the
local chart $(U, \phi)$ we have\footnote{We borrow this observation from
  \cite{MSS} who use it to prove the uniqueness of the extension of the ordinary trace
  on trace-class operators to non integer order operators.}
$$\left({\rm Op}(\partial_{\xi_j} \sigma)u\right)(x)=  \int_{\R^n} e^{i\langle x, \xi\rangle} \partial_{\xi_j}\sigma(x, \xi)
\hat u(\xi) d\,\xi =-i \left(\, {\rm ad}_{x_j} {\rm
Op}(\sigma)u\right)(x)\quad \forall u\in C_{\rm cpt}^\infty (U).$$
Furthermore,
\begin{eqnarray*}
 \rho_\phi\circ \partial_{\xi_j}&=&  L\circ\phi^* \circ {\rm
  Op}\circ \partial_{\xi_j}\\
&=&-i\,   L\circ\phi^*\circ  {\rm ad}_{x_j} \circ {\rm
  Op}\\
&=&-i\,   L\circ  {\rm ad}_{x_j}\circ\phi^* \circ {\rm
  Op}.
\end{eqnarray*}
Since $L$ vanishes on brackets,  $\rho_\phi$ vanishes on derivatives
$\partial_{\xi_j}\tau$.
Similarly, for any $u\in C^\infty_{\rm cpt}(U)$
\begin{eqnarray*}
\left({\rm Op}(\partial_{x_j} \sigma)u\right)(x)&=  &\int_{\R^n}
e^{i\langle x, \xi\rangle} \partial_{x_j}\sigma(x, \xi)
\hat u(\xi) d\,\xi\\
& =&\partial_{x_j} \int_{\R^n} e^{i\langle x, \xi\rangle} \sigma(x, \xi)
\hat u(\xi) d\,\xi -i\int_{\R^n}\xi_j  e^{i\langle x, \xi\rangle} \,\sigma(x, \xi)
\hat u(\xi) d\,\xi\\
&=&\partial_{x_j} \int_{\R^n} e^{i\langle x, \xi\rangle} \sigma(x, \xi)
\hat u(\xi) d\,\xi -\int_{\R^n}  e^{i\langle x, \xi\rangle} \,\sigma(x, \xi)
\widehat{ \partial_{x_j} u(\xi)} d\,\xi\\
&=&\left[\partial_{x_j}, {\rm Op}(\sigma)\right] u(x).
\end{eqnarray*}
Furthermore,   \begin{eqnarray*}
\rho_\phi\circ \partial_{x_j}&=&  L\circ\phi^* \circ {\rm
  Op}\circ \partial_{x_j}\\
&=&   L\circ\phi^*\circ  {\rm ad}_{\partial_{x_j}} \circ {\rm
  Op}\\
&=&  L\circ  {\rm ad}_{\partial_{x_j}}\circ\phi^* \circ {\rm
  Op}.
\end{eqnarray*}
Since $L$ vanishes on brackets it follows that $\rho_\phi$ vanishes on
derivatives $\partial_{x_i} \tau$ and therefore satisfies Stokes' property. \\
By Proposition \ref{prop:uniquenessresU}, $\rho_\phi$ which is
continuous on each $CS_{\rm cpt}^a(\phi(U))$ as a result of the
continuity of $L$ on each $\Cl^a(M)$,  is therefore proportional to the
noncommutative residue so that there is a constant $c_\phi$ such
that
$$\rho_\phi(\sigma)=L(\phi^*{\rm Op}(\sigma))=c_\phi\cdot{\rm res}(\sigma)\quad\forall \sigma \in
CS_{\rm cpt}(\phi(U)).$$
We can now use a partition of unity    $(U_i, \chi_i)_{i \in I}$ subordinated to an atlas $(U_i,
\phi_i)_{i\in I} $ on $M$ to  write any operator  $P\in \Cl(M)$ as a finite sum of localised
operators
$P=\sum_{i\in I} P_i$ with  $ P_i:= \chi_i\, P\,\chi_i.   $ We can further
assume that  $P_i=\phi_i^*{\rm Op}(p_i)$ with $p_i\in CS_{\rm
  cpt}(\phi_i(U_i))$. It follows from the first part of the proof that
 $L(P_i)=\rho_{\phi_i}(p_i) = c_{\phi_i}\cdot {\rm res}(p_i) $ so that
by linearity of $L$, we have
$L(P)=  \sum_{i\in I} L(P_i)= \sum_{i\in I} c_{\phi_i}\cdot {\rm res}(p_i). $
But since the l.h.s is globally defined, the r.h.s is independent of the local
chart; it follows that $L(P)= c\cdot {\rm res}(P)$ for some constant $c\in
\C$.\endsquare
\subsection{Uniqueness of the canonical trace }
\begin{prop}\label{prop:uniquenesscutoffintU} Let ${\cal D}(U)$ be a subset of
$ CS_{\rm cpt}(U)$  containing smoothing symbols which  is stable under multiplication  by smooth
functions:
$$C^\infty_{\rm cpt} (U)\cdot {\cal D}(U)\subset {\cal D}(U).$$
Let
$${\cal S}:=\{\sigma\in  CS_{\rm c.c}(\R^n), \quad f\cdot\sigma\in  {\cal
  D}(U)\quad \forall f\in C^\infty_{\rm cpt}(U)\}.$$
We further assume that $ C^\infty_{\rm cpt} (U)\otimes \left({\cal
  S}\cap CS_{\rm c.c}^a(\R^n)\right)$ is dense in ${\cal D}(U)\cap CS_{\rm
cpt}^a(U)$ for any $a\in\C$  and that it fulfills the
requirements of Theorem \ref{thm:uniquenesscutoffint}.  \\ \\
Then, any continuous linear form \footnote{i.e. its restriction to  ${\cal    D}(U)\cap CS_{\rm cpt}^a(U)$ is continuous for any $a\in\C$.} $\rho:
{\cal D}(U)\to \C$  which satisfies Stokes' property  is proportional to the cut-off regularised integral:
$$\exists c\in \C, \quad \rho(\sigma)=c\cdot \cutoffint_{T^*U}\sigma\quad \forall \sigma\in {\cal D}(U).$$
\end{prop}
\begin{rk}\label{rk:vanishingresdensity}
Since ${\cal S}$ fulfills the requirement of  Theorem
  \ref{thm:uniquenesscutoffint}, we have ${\cal S}\subset {\rm
    Ker}({\rm res})$ and hence ${\rm res}(f\cdot\sigma)=0\quad\forall f\in
  C^\infty_{\rm cpt}(U), \quad \forall\sigma \in  {\cal S}$. By a density argument using the continuity of the
  residue on symbols of constant order, this implies that ${\rm res}(f\cdot
  \sigma)=0 \quad\forall f\in
  C^\infty_{\rm cpt}(U), \quad \forall \sigma\in {\cal D}(U).$ The
  requirements of the proposition therefore imply that symbols in ${\cal D}(U)$ have
  vanishing residue density ${\rm
    res}_x(\sigma(x, \cdot))=0\quad \forall \sigma\in {\cal D}(U)\quad\forall x\in U.$
\end{rk}
{\bf Proof:} We closely follow the proof of Proposition \ref{prop:uniquenessresU}.\\
For a fixed $f\in C^\infty_{\rm cpt}(U)$ the map 
$\sigma\mapsto  \rho(f\, \sigma)$ defines a continuous linear form on
${\cal S}$
which vanishes on derivatives in $\xi$ since we have $\rho(f\, \partial_{\xi_j} \sigma)=  \rho(\partial_{\xi_j}(f\,  \sigma))=0$
for any smooth function $f\in C^\infty_{\rm cpt}(U)$.
By Theorem \ref{thm:uniquenesscutoffint}, it follows that  there is a constant $c(f)$
such that $\rho(f\, \sigma)= c(f)\cdot \cutoffint_{\R^n}\sigma$ for any $\sigma\in
 {\cal S}$ and for any  $
 f\in C_{\rm cpt}^\infty(U).$ Since
$CS_{\rm ct}^a(U)\otimes\left( {\cal
  S}\cap CS_{\rm c.c}^a(\R^n)\right)$  is dense in
${\cal D}(U)\cap CS_{\rm cpt}^a(U)$ for any $a\in \C$ and since $\rho$ is continuous
on ${\cal D}(U)\cap CS_{\rm cpt}^a(U)$ it  follows that 
$$\rho(\sigma)= F\left( \cutoffint_{\R^n}\sigma(x, \cdot)\right)$$ for some
continuous distribution $F:f\mapsto \int_U f(x)\, \phi(x)\, dx $ with $\phi \in  C^\infty(U)$.
From the closedness of $\rho$ we infer that $ \rho(\partial_{x_i}f\,
\sigma)=F\left(\partial_{x_i}f\, 
  \cutoffint_{\R^n} \sigma\right)=0$ for any $\sigma \in {\cal S}$ and any $f\in C_{\rm
  cpt}^\infty(U)$. Choosing $\sigma$ such that  $\cutoffint_{\R^n} \sigma\neq
0$ implies that $F(\partial_{x_i}f)=0$ and hence that  $\phi$ is  constant and 
$$\rho(\sigma)=c\cdot  \int_U \cutoffint_{\R^n}\sigma=  c\cdot \cutoffint_{U\times \R^n}\sigma\quad \forall 
\sigma \in {\cal S}(U).$$ \endsquare
\begin{ex} ${\cal D}(U):= CS_{\rm cpt}^{\notin\Z}(U)$ satisfies the
  assumptions of the proposition. Indeed in that case ${\cal S}= CS_{\rm
    c.c}^{\notin\Z}(\R^n)$ fulfills the
  requirements of Theorem \ref{thm:uniquenesscutoffint} and   $C^\infty_{\rm
    cpt}(U)\otimes  CS_{\rm c.c}^a(\R^n)$ is
  dense in $  CS_{\rm cpt}^a(U)$ for any non
  integer order $a$.
\end{ex}
\begin{ex}If the dimension  $n$ is odd  then $${\cal D}(U):= CS_{\rm cpt}^{odd}(U)=
  \{\sigma\in CS_{\rm cpt}(U), \quad\sigma_{a-j}(x, -\xi)=
  (-1)^{a-j}\,\sigma_{a-j}(x, \xi)\quad\forall x\in U\quad {\rm with}\quad a={\rm ord}\sigma\}$$ satisfies the
  assumptions of the proposition. Indeed, in that case ${\cal S}= \{\sigma\in CS_{\rm c.c}(\R^n), \quad\sigma_{a-j}( -\xi)=
  (-1)^{a-j}\,\sigma( \xi) \quad {\rm with}\quad a={\rm ord}\sigma\} $  fulfills the
  requirements of Theorem \ref{thm:uniquenesscutoffint} and  $C^\infty_{\rm
    cpt}(U)\otimes \left({\cal S}\cap CS_{\rm c.c}^a(\R^n)\right)$ is
  dense in $ CS_{\rm cpt}^{odd}(U)\cap  CS_{\rm cpt}^a(U)$ for any $a\in \C$.
\end{ex}
\begin{ex}A similar statement holds in the even dimensional case replacing odd
  class symbols by even class symbols $\sigma\in CS_{\rm cpt}^{even}(U)$ i.e by
  symbols $\sigma$  that satisfy  the requirement $\sigma_{a-j}(x, -\xi)=
  (-1)^{a-j+1}\,\sigma_{a-j}(x, \xi)$ for any $ x\in U$.
\end{ex}

\begin{thm}\label{thm:uniquenesstraceM} Let  ${\cal D}(M)$ be
 a subset of $\Cl(M)$ containing smoothing operators which is  stable under multiplication by smooth functions:
$$\Ci(M)\cdot {\cal D}(M)\subset {\cal D}(M).$$ We further assume that
\begin{enumerate}
\item the canonical
trace is well-defined on ${\cal D}(M)$ and vanishes on brackets in ${\cal
  D}(M)$,
\item  given any local
chart $(U, \phi)$ on $M$
$${\cal S}_\phi:=\{\sigma\in  CS_{\rm c.c}(\R^n), \quad \phi^*{\rm
  Op}(f\cdot\sigma)\in {\cal
  D}(M)\quad \forall f\in C^\infty_{\rm cpt}(\phi(U))\}$$ fulfills the
assumptions of Theorem \ref{thm:uniquenesscutoffint}.
\end{enumerate}
 Then any  continuous \footnote{i.e. which restricts to a  continuous map   on ${\cal }D(M)\cap \Cl^a(M)$ for any $a\in \C$.}   linear form \footnote{By linear we mean here that $L(\alpha A+ \beta B)=
\alpha  L( A)+ \beta L(B)$ whenever $A, B,\alpha A+ \beta B\in  {\cal D}(M)$ }
$L:{\cal D}(M)\to \C$
  which  vanishes on brackets:
\begin{equation}\label{eq:vanishesonbrackets}L\left([A, B]\right)=0\quad \forall A, B\in \Cl(M)\quad {\rm s.t.}\quad
[A,B]\in {\cal D}(M) \end{equation} is proportional to the canonical trace:
$$\exists c\in\C, \quad L(A)= c\cdot {\rm TR}(A)\quad \forall A\in {\cal D}(M).$$
\end{thm}  \begin{rk} According to Remark \ref{rk:vanishingresdensity}, it  follows from the  assumption 2
  that ${\cal D}(M)$ is contained in
$${\rm Ker}({\rm res})_{\rm loc}(M):= \{A\in \Cl(M),\quad {\rm s.t.} \quad  f\, \cdot A\in
  {\rm Ker}({\rm res})\quad \forall f\in \Ci(M)\}$$ which corresponds to the   linear space of operators $A\in \Cl(M)$ with vanishing residue density
  i.e.:
$${\rm Ker}({\rm res})_{\rm loc}(M)= \{A\in \Cl(M),\quad {\rm s.t.}\quad \sigma_A(x, \cdot)\in
  {\rm Ker}({\rm res}_x)\quad \forall x\in M\}$$
since
\begin{eqnarray*}
&{}& {\rm res}(f\, A)=0\quad \forall f\in \Ci(M)\\
&\Leftrightarrow& \int_M f(x)\, {\rm res}_x\left(\sigma_A\right)(x, \cdot)\,
dx=0\quad \forall f\in \Ci(M)\\
&\Leftrightarrow & {\rm res}_x\left(\sigma_A\right)(x, \cdot)=0\quad \forall
x\in M.
\end{eqnarray*}
\end{rk}{\bf Proof:} Since the proof closely follows that of
Theorem \ref{thm:uniquenessresM}, we do not repeat some of the steps  common
to the two proofs.\\
Let us first observe that given any local chart $(U, \phi)$ on $M$   the set  $${\cal D}(\phi(U)):= \{\sigma \in
CS_{\rm cpt}(\phi(U)), \quad  \phi^*\circ {\rm
  Op}(\sigma)\in {\cal D}(M)\}$$  fulfills the assumptions of Proposition
\ref{prop:uniquenesscutoffintU} with $U$ replaced by $\phi(U)$ and with ${\cal
  S}$ replaced by
${\cal S}_\phi$ as in the statement of the theorem. \\
From a linear form $L$ on ${\cal D}(M)$ which obeys the requirements of the
theorem  we can build the   linear form  $$\rho_\phi:=L\circ\phi^*\circ {\rm
  Op}$$  on ${\cal D}(\phi(U))$ which obeys the requirements of Proposition  \ref{prop:uniquenessresU}.
Hence $\rho_\phi$  is proportional to the cut-off regularised integral  so that there is a
constant $c_\phi$ such that
$$\rho_\phi(\sigma)=L(\phi^*{\rm
  Op}(\sigma))=c_\phi\cdot\cutoffint_{T^*\phi(U)}\sigma\quad\forall \sigma \in
{\cal D}(\phi(U)).$$
As before, using a partition of unity to  write any operator  $P\in \Cl(M)$ as a finite sum of localised
operators
$P=\sum_{i\in I} P_i$ with  $ P_i:= \chi_i\, P\,\chi_i,   $  with   $P_i=\phi_i^*{\rm Op}(p_i)$ with $p_i\in CS_{\rm
  cpt}(\phi_i(U_i))$ we infer that
 $L(P_i)=\rho_{\phi_i}(p_i) = c_{\phi_i}\cdot \cutoffint_{T^*\phi_i(U_i)} p_i $ so that
by linearity of $L$
$$L(P)=  \sum_{i\in I} L(P_i)= \sum_{i\in I}
c_{\phi_i}\cdot \cutoffint_{T^*\phi_i(U_i)} p_i  . $$
But since the l.h.s is globally defined, the r.h.s is independent of the local
chart; it follows that $L(P)= c\cdot \cutoffint_{T^*M}\sigma(P)=c\cdot {\rm TR}(P)$ for some constant $c\in
\C$.\endsquare\\  \\
Here are a few known examples of sets ${\cal D}(M)$ which obey assumptions 1
and 2 of the above theorem. In particular, they lie in  ${\rm Ker(res)}_{\rm loc}(M)$.
\begin{ex}The set   $\Cl^{\notin \Z}(M)$ of non integer order classical pseudodifferential
  operators on $M$.
\end{ex}
\begin{ex} The set $$\Cl^{\rm odd}(M)= \{ A\in \Cl^\Z(M), \quad
\sigma(A)\sim \sum_{j=0}^\infty \chi \, \sigma_{a-j}(A), \quad
\sigma_{a-j}(A)(x, -\xi)=(-1)^{a-j}\sigma_{a-j}(A)(x, \xi) \}$$ of  odd-class operators
on odd
dimensional manifolds $M$ introduced in \cite{KV} (see also \cite{Gr} where
such  operators are called even-even).
\end{ex}
\begin{ex}The set $$\Cl^{\rm even}(M)= \{ A\in \Cl^\Z(M), \quad
\sigma(A)\sim \sum_{j=0}^\infty \chi \, \sigma_{a-j}(A), \quad
\sigma_{a-j}(A)(x, -\xi)=(-1)^{a-j+1}\sigma_{a-j}(A)(x, \xi) \}$$ of   even-class operators
on even
dimensional manifolds $M$ (see \cite{Gr} where
such  operators are called even-odd).
\end{ex}
 Applying Theorem \ref{thm:uniquenesstraceM} to ${\cal D}(M)=  \Cl^{\notin
    \Z}(M)$, resp. ${\cal D}(M)=\Cl^{\rm odd}(M)$ in odd dimensions,  resp. ${\cal D}(M)=\Cl^{\rm
    even}(M)$ in even dimensions, leads to the following uniqueness result.
\begin{cor} The canonical trace is (up to a multiplicative constant) the unique linear form on $\Cl^{\notin
    \Z}(M)$, resp. $\Cl^{\rm odd}(M)$ in odd dimensions,  resp. $\Cl^{\rm
    even}(M)$ in even dimensions which is continuous on operators of constant
  order and which vanishes on brackets that lie in $\Cl^{\notin
    \Z}(M)$, resp. $\Cl^{\rm odd}(M)$ in odd dimensions,  resp. $\Cl^{\rm
    even}(M)$.
\end{cor}
\begin{rk} In the course of the proof we showed that the vanishing of $L$ on brackets
(\ref{eq:vanishesonbrackets}) implies Stokes' property for $
\rho_\phi$. Conversely, Stokes' property for  $
\rho_\phi$ implies that $L(\phi^*{\rm Op}(\sigma)):= \rho_\phi(\sigma)$ vanishes on brackets
$[x_i, \cdot]$ and $[\partial_{x_i}, \cdot]$ contained in ${\cal D}(M)$. But  this implies that $L$ vanishes on
brackets $[P_U,\cdot]\in {\cal D}(M)$ where $P_U$ is the localisation of any
classical pseudodifferential operator. Stokes' property on
symbols and the vanishing on brackets of operators are therefore equivalent.
\end{rk}
\vfill \eject 
\end{document}